\documentclass[a4paper,10pt]{article}

\usepackage{latexsym}
\usepackage{amsmath,amsfonts,amssymb}
\usepackage{graphicx,color,epsfig}

\def\C{\mathbb{C}}
\def\N{\mathbb{N}}

\def\R{\mathbb{R}}
\def\T{\mathbb{T}}
\def\Z{\mathbb{Z}}

\def\leq{\leqslant}
\def\geq{\geqslant}
\def\eps{\epsilon}
\newcommand{\pa}{\partial}

\newtheorem{Def}{Definition}[section]
\newtheorem{Thm}[Def]{Theorem}
\newtheorem{Lem}[Def]{Lemma}
\newtheorem{Pro}[Def]{Proposition}
\newtheorem{Cor}[Def]{Corollary}
\newtheorem{Exam}[Def]{Example}

\begin{document}

\title{Landau damping to partially locked states in the Kuramoto model}
\author{Helge Dietert$^1$, Bastien Fernandez$^{2}$\ and David G\'erard-Varet$^3$}
\date{June 14, 2016}
\maketitle
\begin{center}
$^1$ Department of Pure Mathematics and Mathematical Statistics\\
University of Cambridge\\
Wilberforce Road\\
Cambridge CB3 0WA UK\\
\medskip
$^2$ Laboratoire de Probabilit\'es et Mod\`eles Al\'eatoires\\
CNRS - Universit\'e Paris 7 Denis Diderot - UPMC\\
75205 Paris CEDEX 13 France\\
\medskip
$^3$ Institut de Math\'ematiques de Jussieu - Paris Rive Gauche\\
Universit\'e Paris 7 Denis Diderot - Sorbonne Paris Cit\'e\\
75205 Paris CEDEX 13 France
\end{center}

\begin{abstract}
  In the Kuramoto model of globally coupled oscillators, partially
  locked states (PLS) are stationary solutions that incorporate the
  emergence of partial synchrony when the interaction strength
  increases. While PLS have long been considered, existing results on
  their stability are limited to neutral stability of the linearized
  dynamics in strong topology, or to specific invariant subspaces
  (obtained via the so-called Ott-Antonsen (OA) ansatz) with specific
  frequency distributions for the oscillators. In the mean field
  limit, the Kuramoto model shows various ingredients of the Landau
  damping mechanism in the Vlasov equation. This analogy has been a
  source of inspiration for stability proofs of regular Kuramoto
  equilibria. Besides, the major mathematical issue with PLS
  asymptotic stability is that these states consist of heterogeneous
  and singular measures.  Here, we establish an explicit criterion for
  their spectral stability and we prove their local asymptotic
  stability in weak topology, for a large class of analytic frequency
  marginals. The proof strongly relies on a suitable functional space
  that contains (Fourier transforms of) singular measures, and for
  which the linearized dynamics is well under control.  For
  illustration, the stability criterion is evaluated in some standard
  examples. We show in particular that no loss of generality results
  in assuming the OA ansatz.  To our best knowledge, our result
  provides the first proof of Landau damping to heterogeneous and
  irregular equilibria, in absence of dissipation.
\end{abstract}

\section{Introduction}
\subsection{Landau damping in the Kuramoto model}
In its original version, the Kuramoto model is a simple (finite-dimensional) model of globally coupled oscillators \cite{ABP-VRS05,S00}. Its dynamics for a population of $N$ oscillators can be regarded as special case of the Winfree model \cite{W67} and is given by the following set of coupled first order equations \cite{K75,K84}
\[
\dot{\theta_i}= \omega_i + \frac{K}{N} \sum_{j=1}^N \sin(\theta_j - \theta_i),\ \forall i = 1,\dots,N,
\]
where $\theta_i\in\T=\R/(2\pi\Z)$ are the oscillator phases, $\omega_i\in \R$ are their frequencies and
$K\geq 0$ parametrizes the strength of the coupling term.

With such simple ingredients, this model has become a paradigm of the transition to synchrony in collective systems, and has been employed in numerous disciplines such as Chemistry, Biology, Social Sciences, etc \cite{ABP-VRS05,PA15}. In spite of substantial attention, the mathematically rigorous description of its asymptotic dynamics mostly remains incomplete, especially when the frequencies are randomly distributed, as in the original formulation. Hence, focus is made on the dynamics of empirical measures $\frac1{N}\sum_i\delta_{\theta_i,\omega_i}$ and their mean-field limit (continuum approximation).

In this setting, the distribution of a (finite or infinite) population at time $t$ is described by a probability measure $f(t)$ on the cylinder $\T\times\R$ \cite{G13}. The measure $f(t)$ is a weak solution of the Kuramoto equation (first suggested by Sakaguchi \cite{S88} in a more general setting, see also \cite{SM91})
\[
\pa_t f + \pa_\theta (fV[f])= 0,
\]
where
\[
V[f](\theta,\omega)=\omega + K \int_{\T\times\R} \sin(\theta' - \theta) f(d\theta',d\omega'),
\]
for $(\theta,\omega)\in\T\times\R$. The Cauchy problem for the weak formulation of the Kuramoto equation is well-posed \cite{L05}. Moreover, the dynamics preserves the frequency marginal, {\sl i.e.}\ we have $\int_\T f(t,d\theta,d\omega)=g(d\omega)$ for all $t\geq 0$ (where for the sake of notation, $g$ is identified with its density in the absolutely continuous case).
Both $K$ and $g$ are crucial parameters of this model. (NB: In principle, $K$ could be absorbed in $g$ via a rescaling of time and frequencies. However, we choose to keep the parameters separated in agreement with the original formulation of the model.)

The Kuramoto model and its variations have rich and diverse phenomenology depending upon $K$ and $g$ \cite{ABP-VRS05,PA15}. Of special interest is the order parameter $r(t)$ defined by
\[
r(t)=\int_{\T\times\R}e^{i\theta}f(t,d\theta,d\omega).
\]
The modulus $|r(t)|$ quantifies coherent behavior in the population; this quantity is small when the population is splayed out, and it is near 1 when the population is close to full synchrony.  Moreover, the Kuramoto equation can be expressed in terms of the order parameter, {\sl viz.}
$$ \pa_t f + \pa_\theta  \left( \omega f + \frac{K}{2i} (e^{- i\theta}r - e^{i \theta} \overline{r})f \right) = 0. $$
In the simplest case of an even and unimodal\footnote{{\sl i.e.}\ unimodal means that $g$ is monotonically increasing up to its maximum and then monotonically decreasing} function $g$, numerical simulations together with heuristic arguments
have identified two distinct phases depending on $K$, which can be
summarized as follows (see \cite{S00} for more details).
\begin{itemize}
\item[$\bullet$] The incoherent (homogeneous) stationary state $f_\text{inc}(\theta,\omega)=\frac{g(\omega)}{2\pi}$ is stable for $K<K_c=\frac{2}{\pi g(0)}$, and asymptotic damping of $r(t)$ results as $t\to+\infty$, for trajectories issued from typical initial conditions.
\item[$\bullet$] The incoherent state undergoes a supercritical pitchfork bifurcation at $K=K_c$ that generates a stable (inhomogeneous) partially locked state (PLS, see definition below). As a consequence, for $K>K_c$, the quantity $|r(t)|$ typically approaches a positive value as $t\to+\infty$, {\sl i.e.}\ coherent behavior emerges in the long term. Moreover, the limit $|r(+\infty)|$ increases with $K$ and tends to 1 when $K\to +\infty$.
\end{itemize}
For more general distributions $g$, the expression of $K_c$ has to be adapted \cite{D14,FG-VG16} and the bifurcation can be of different type \cite{MBSOSA09}.

While this phenomenology had been identified in the early studies
\cite{K75,K84}, its mathematical proof has resisted achievement until
very recently \cite{C15,D14,FG-VG16}, and yet the results are limited
to the asymptotic stability of the incoherent state and to the
neighborhood of the bifurcation at $K=K_c$. (Independently, in the
case where $g$ is concentrated on a single frequency, dissipation implies
asymptotic stability of the fully synchronized state $r=1$ for every $K>0$, both in the finite dimensional Kuramoto
model \cite{CHJK12} and its continuum limit \cite{CCHKJ14}).

The main technical issue comes from the transport term $\omega\partial_\theta f$ in the equation. This term implies that in strong topology, the operator associated with the linearized perturbation dynamics has continuous spectrum on the imaginary axis \cite{MS07}; hence the standard proof of asymptotic stability cannot apply in this setting. In fact, when evaluated using $L^2$(-like) norms, perturbation sizes increase with time \cite{D14,N15}!

A continuous spectrum on the imaginary axis is reminiscent of the Vlasov equation where the same transport term occurs. The similarity with that equation is actually stronger; in both cases, any homogeneous density (or measure), {\sl i.e.}\ not depending on the angular/spatial variable, is a stationary solution. Also, the integrated linearized dynamics of the Fourier modes of the corresponding perturbations are both given by a Volterra equation of the second kind. In addition, the asymptotic relaxation of the order parameter $r(t)$ is the analogue of the Landau damping in the Vlasov setting \cite{SMM92}.

Accordingly, the proofs in \cite{D14,FG-VG16} were inspired by the
proof of Landau damping, which was first fully completed in
\cite{MV11} for the Vlasov-Poisson equation (see also \cite{V10} for a
didactic presentation). In particular, the proof in \cite{FG-VG16}
closely follows the bootstrap argument developed for the Vlasov-HMF
model in \cite{FR16} (NB: The potential in the Vlasov-HMF model
consists of the first Fourier mode, as in the Kuramoto model, instead of being singular as in the
Vlasov-Poisson equation).

In spite of these results, to the best of our knowledge, a rigorous
portrayal of the dynamics in the neighborhood of a Kuramoto stationary
state for arbitrary $K>K_c$ is still missing. The aim of this paper is
to fill this gap and, to a broader extent, to prove asymptotic
stability of PLS under an appropriate spectral
condition.\footnote{Besides, the Kuramoto equation is indifferent to frequency
  translations, {\sl i.e.}\ if $f(t,\theta,\omega)$ is a solution,
  then, for every $\Omega\in\R$, $f(t,\theta+\Omega t,\omega+\Omega)$ is a also solution, with marginal density
  $g(\omega+\Omega)$. Therefore, up to a translation of $g$, the
  existence and stability of any purely rotating solution
  $t\to f(\theta+\Omega t,\omega)$ is equivalent to the existence and
  stability of the stationary state $f$. In other words, our results
  also apply to the dynamics of rotating trajectories and their
  perturbations.}

\subsection{The partially locked states}\label{S-PLS}
A partially locked state (PLS) is a stationary solution $f_\text{pls}$ of the Kuramoto equation with non-zero order parameter $r=r(t)\neq 0$. The equation commutes with any rotation of angle $\Theta$ on the circle, or more precisely, with their representation $R_\Theta$ on measures on the cylinder.
These transformations only affect the phase of the order parameter, {\sl viz.}\ $\int_{\T\times\R}e^{i\theta}R_\Theta f(d\theta,d\omega)=e^{-i\Theta}\int_{\T\times\R}e^{i\theta}f(d\theta,d\omega)$. Accordingly, PLS come in one-parameter families, namely circles of the form $\{R_\Theta f_\text{pls}\}_{\Theta\in\T}$, where we can assume w.l.o.g.\
\begin{equation}
r_\text{pls}:=\int_{\T\times\R}e^{i\theta}f_\text{pls}(d\theta,d\omega)\in\R^+.
\label{ORDPARAM}
\end{equation}
Using this assumption, we can solve the stationary state equation
$$ \pa_\theta ( (\omega  -  K  r_\text{pls}  \sin \theta) f ) = 0 $$
 (for a detailed reasoning, see {\sl e.g.}\ \cite{ABP-VRS05,MS07,S00}) and get the following expression
\begin{equation}
f_\text{pls}(\theta,\omega)=\left\{\begin{array}{ccl}
\left(\alpha(\omega)\delta_{\arcsin(\frac{\omega}{Kr_\text{pls}})}(\theta)+(1-\alpha(\omega))\delta_{\pi-\arcsin(\frac{\omega}{Kr_\text{pls}})}(\theta)\right)g(\omega)&\text{if}&|\omega|\leq Kr_\text{pls}\\
\frac{\sqrt{\omega^2 - (Kr_\text{pls})^2}}{2\pi |\omega - Kr_\text{pls}\sin \theta|} g(\omega) &\text{if}&|\omega|> Kr_\text{pls}.
\end{array}\right.
\label{EXPREPLS}
\end{equation}
Here the measurable function $\alpha:[-Kr_\text{pls},Kr_\text{pls}]\to [0,1]$ is arbitrary and describes the combination of the masses associated with the two equlibria $\arcsin(\frac{\omega}{Kr_\text{pls}})$ and $\pi-\arcsin(\frac{\omega}{Kr_\text{pls}})$, of the equation $\dot\theta=\omega-Kr_\text{pls}\sin\theta$. This equation governs the oscillator dynamics on $\T$ in the original model. Obviously, these equilibria only exist for $|\omega|\leq Kr_\text{pls}$, and when in these states, the oscillators are ``locked" to their frequency. On the opposite, for $|\omega|> Kr_\text{pls}$, the equation does not have any equilibrium point  and the oscillators must rotate forever; hence the absolutely continuous distribution in this range.

Imposing that the solution \eqref{EXPREPLS} satisfies assumption \eqref{ORDPARAM} yields the following self-consistency equation for the order parameter $r_\text{pls}$ (see \cite{MS07} or use Fourier coefficients, expression \eqref{FOURCOEF} in Appendix \ref{A-STAB})
\[
K\int_{-1}^1(2\alpha(Kr_\text{pls} \omega)-1)\sqrt{1-\omega^2}g(Kr_\text{pls}\omega)d\omega=1,
\]
and also the condition
\[
\int_{1}^{+\infty}\sqrt{\omega^2-1}\left(g(Kr_\text{pls}\omega)-g(-Kr_\text{pls}\omega)\right)d\omega=\int_\R\omega g(Kr_\text{pls}\omega)d\omega,
\]
which obviously holds when $g$ is even. These constraints play the role of existence conditions for PLS.

Let $f_\text{s}$ denote a PLS with $\alpha=1\ \text{a.e.}$ and let $r_\text{s}>0$ be the corresponding order parameter. A basic continuity argument on $r_\text{s}$ across $[0,1]$ in the consistency equation shows that, for any continuous and even density $g$, a PLS $f_\text{s}$ exists for every $K>K_c$. (If $g$ is unimodal, then
monotonicity implies uniqueness and that $K>K_c$ is necessary for
existence. If otherwise, 0 is a local minimum for $g$, then two (or
more) PLS $f_\text{s}$ exist for $K$ in a left-neighborhood of
$K_c$).

When it comes to stability, the analysis of the equation $\dot\theta=\omega-Kr_\text{pls}\sin\theta$, suggests that only PLS of the type $f_\text{s}$ can be stable \cite{MS07}. As shown in Appendix \ref{A-STAB}, this choice of all locked oscillators at the stable stationary point $\arcsin(\frac{\omega}{Kr_\text{pls}})$ can actually be justified as a PLS stability requirement (see Proposition \ref{UNIQPLS}).

As for the stability analysis itself, studies can be found in the literature on spectral properties of the linearized perturbation dynamics, either in the full space when $g$ is symmetric and unimodal \cite{MS07,SM91}, or in the so-called Ott-Antonsen manifold\footnote{The Ott-Antonsen manifold is the set of those probability measures on the cylinder which write
$$
\frac{g(\omega)}{2\pi} \sum_{\ell=0}^\infty \left( h(\omega)^\ell  e^{i \ell \theta} + \overline{h(\omega)}^\ell e^{-i \ell \theta} \right),
$$
for some amplitude function $h$, see Section \ref{S-OA} for details. Remarkably, this set is invariant under the Kuramoto flow.}
when $g$ is arbitrary with compact support \cite{OW13}. All these works consider strong topology. As previously mentioned, they can only result in neutral stability and this prevents any control of the nonlinear terms. (NB: When $g$ is a rational function, the evolution of the  order parameter for solutions in the Ott-Antosen manifold can be described through a  finite-dimensional system; hence no such issue exists and convergence of the order parameter follows from spectral stability \cite{MBSOSA09}.)
Besides, transport terms are known to produce a stabilizing effect on the linearized dynamics, when evaluated in weak topology \cite{C15,D14,MV11}. This suggests to consider that topology instead of the strong one.

Accordingly, the main lines of our PLS asymptotic stability proof are
as follows. We first introduce a weak functional space (after passing
to Fourier variables) that turns out particularly well-suited. Indeed,
we shall not only prove that this space contains the PLS $f_s$ and the initial value problem is globally well-posed therein, but also that the essential spectrum of the (complexified) linearized operator around $f_s$ entirely lies in the left part of the complex plane. Then, we establish a criterion that controls the remaining discrete spectrum, namely we ensure that 0 is the only eigenvalue with non-negative real part (and is simple).

That 0 is always an eigenvalue of the linearized operator is a consequence of the rotation symmetry. This suggests that asymptotic stability should apply to the whole PLS circle $\{R_\Theta f_\text{s}\}_{\Theta\in\T}$, {\sl i.e.}\ that any sufficiently small perturbation of any state $R_\Theta f_\text{s}$ asymptotically relaxes, under the full nonlinear dynamics, to $R_{\Theta'} f_\text{s}$ for some $\Theta'\in\T$ close to $\Theta$. This is exactly what we prove to happen when the stability criterion holds. The proof proceeds by projecting out the coordinate along the PLS circle, and show asymptotic stability of the resulting relative equilibrium, see {\sl e.g.}\ \cite{CL00,HI11}.

The analysis holds for a large class of analytic frequency distributions $g$. In addition, various examples of the Kuramoto literature are re-visited in Section \ref{S-STAB}, which also includes stability considerations in the Ott-Antonsen manifold.

As intended, our result in particular completes the mathematical description of the Kuramoto equation for symmetric and unimodal marginal densities, and for arbitrary interaction strength. Moreover, to our best knowledge, it provides the first
spectral stability based proof of asymptotic convergence (in the weak sense) to inhomogeneous and irregular stationary states in the mean field limit of a system of interacting particles. We hope that in the future, this approach can be extended to other systems and in particular to the asymptotic stability of inhomogeneous states in the Vlasov-HMF model, as reported in {\sl e.g.}\ \cite{BOY15}.

\section{Main result}
As announced above, the core analysis in this paper operates on Fourier transforms of measures on the cylinder defined by
\[
\widehat{f}_\ell(\tau)=\int_{\T\times\R}e^{-i(\ell\theta+\tau\omega)}f(d\theta,d\omega),\ \forall (\ell,\tau)\in\Z\times\R.
\]
That the Kuramoto dynamics preserves frequency marginals implies $\widehat{f}_0(t,\tau)=\widehat{g}(\tau)$ for all $t\geq 0$ for the Fourier transforms (where, evidently, $\widehat{g}(\tau)=\int_{\R}e^{-i\tau\omega}f(d\omega)$).
Therefore, given also that the solutions $f(t)$ are real measures, in order to get information on the whole $\widehat{f(t)}$, all we need to control is the restriction $\widehat{f(t)}|_{\N\times\R}$. For the sake of notations, we shall use the symbols $\widehat{f}$ (resp.\ $\widehat{f(t)}$) to denote $\widehat{f}|_{\N\times\R}$ (resp.\ $\widehat{f(t)}|_{\N\times\R}$), throughout the paper.

Inspired by the approach to the incoherent state stability in \cite{D14}, the following weighted norms will appear convenient. Let $a > 0$ and $k \in \R$ be arbitrary.  Given $h : \R \to \C$, we consider
\begin{equation}
\|h\|_{a}=\left(\int_\R e^{2a \tau}\left(|h(\tau)|^2+|h'(\tau)|^2\right)d\tau\right)^{\frac12}.
\label{SOBNORM}
\end{equation}
Moreover, given $u:\N\times\R\to\C$, we consider
\[
\|u\|_{a,k}=\left(\sum_{\ell\in\N}\int_\R e^{2a \tau}\ell^{2k}\left(|u_\ell(\tau)|^2+|\partial_\tau u_\ell(\tau)|^2\right)d\tau\right)^{\frac12}.
\]
The weight $e^{2a\tau}$ here, as opposed to the standard choice $e^{2a|\tau|}$, allows to include Fourier transforms of irregular states, especially $\widehat{f_s}$, see Appendix \ref{A-STAB}. Also, convergence in the  norm $\|\cdot\|_{a,0}$ implies weak convergence of probability measures with fixed frequency marginals, see Appendix \ref{A-WEAKCONV}. Moreover and as indicated above, using a $L^2$-norm with exponential weight results in a shift to the stable half-space, of the essential spectrum of the linearized generator at $\widehat{f_s}$. Of note, a similar idea has already been applied to the stability of KdV solitons \cite{PW94}.

Obviously, a prerequisite is to make sure that the Cauchy problem for the Kuramoto equation is well-posed for data whose Fourier transform have finite $\|\cdot\|_{a,0}$ norm, and   that we have $\|\widehat{f_\text{s}}\|_{a,0}<+\infty$ for the PLS with $\alpha=1\ \text{a.e.}$ This is addressed in Section \ref{S-WELLPOSED} and Appendix \ref{A-STAB} respectively.

Now, in order to state the stability condition, we need to introduce
some notations. Given $\lambda\in\C$ with $\text{Re}(\lambda)\geq 0$
and $r\in\R^+$, let $M(\lambda,r)$ be the
$2\times 2$ matrix defined by
\[
M(\lambda,r)=\left(\begin{array}{cc}
J_0(\lambda,r)&J_2(\lambda,r)\\
\overline{J_2(\bar\lambda,r)}&\overline{J_0(\bar\lambda,r)}\end{array}\right),
\]
where\footnote{Strictly speaking, the integrals here are only well-defined for $\text{Re}(\lambda)>0$. For $\text{Re}(\lambda)=0$, the quantities $J_k$ are defined by using continuity, see Lemma \ref{EQUIVM} below.}
\begin{equation}
J_k(\lambda,r)=\int_\R\frac{\beta^k\left(\frac{\omega}{Kr}\right)}{\lambda+ i\omega+Kr\beta\left(\frac{\omega}{Kr}\right)}g(\omega)d\omega,\ \text{for}\ k\in\N\cup\{0\},
\label{DEFJK}
\end{equation}
and $\beta$ is defined on $\R$ by
\begin{equation}
\beta(\omega)=-i\omega+\left\{\begin{array}{ccl}
\sqrt{1-\omega^2}&\text{if}&|\omega|\leq 1\\
i\omega\sqrt{1-\omega^{-2}}&\text{if}&|\omega|> 1.
\end{array}\right.
\label{BETA}
\end{equation}
\begin{Thm}\label{thmNL}
Assume that $g\in L^\infty(\R)$, $\|\widehat{g}\|_a <+\infty$ for some $a>0$ and a PLS $f_\text{s}$ with marginal density $g$ and order parameter $r_\text{s}\in\R^+$ exists and satisfies the following conditions:
\begin{equation}
\left\{
\begin{aligned}
& \det \left(\text{\rm Id}-\frac{K}2 M(\lambda,r_\text{s})\right)\neq 0,\ \forall \lambda\neq 0\ \text{with}\ \text{\rm Re}(\lambda)\geq 0, \\
&  \:   \liminf_{\lambda \rightarrow 0}  \left| \frac{1}{\lambda} \det \left(\text{\rm Id}-\frac{K}2 M(\lambda,r_\text{s})\right)\right| > 0.
\end{aligned}
\right.
\label{STABCOND}
\end{equation}
Then, there exists $\epsilon,b>0$ such that for every probability measure $f_\text{\rm in}$ with marginal density $g$ and satisfying
\[
\|\widehat{f_\text{\rm in}}-\widehat{f_\text{s}}\|_{a,0}<\epsilon,
\]
there exists $\Theta_\infty\in\T$ so that the solution $t\mapsto f(t)$ of the Kuramoto equation with initial data $f(0)=f_\text{\rm in}$ has the following asymptotic behavior
\[
\|\widehat{f(t)}-\widehat{R_{\Theta_\infty}f_\text{s}}\|_{a,0}=O(e^{-bt}).
\]
In particular, we have
\[
\lim_{t\to+\infty}f(t)=R_{\Theta_\infty} f_\text{s},
\]
in weak topology.
\end{Thm}
(Theorem \ref{thmNL} is proved in Section \ref{S-NLSTAB}.) 
Asymptotic stability of the whole circle $\{R_\Theta f_\text{s}\}_{\Theta\in\T}$ immediately follows from commutation with the rotations $R_\theta$. Moreover, notice that global stability of PLS does not hold because the incoherent state $f_\text{inc}(\theta,\omega)=\frac{g(\omega)}{2\pi}$ always exists as a stationary state and we have $\|\widehat{f_\text{inc}}-\widehat{f_\text{s}}\|_{a,0}<+\infty$. Hence, the size of PLS perturbations always has to be limited.

We shall see in Section \ref{S-LINEAR} that the linearized generator at $\widehat{f_\text{s}}$ must have finitely many eigenvalues in the half-space $\text{Re}(\lambda)>-a$ and that these eigenvalues are determined by the equation $\det \left(\text{\rm Id}-\frac{K}2 M(\lambda,r_\text{s})\right)=0$. As already mentioned, symmetry considerations imply that 0 must be such an eigenvalue.\footnote{One can also directly check that $\det \left(\text{\rm Id}-\frac{K}2 M(0,r_\text{s})\right)=0$, see beginning of Section \ref{S-STAB}.} In other words, the first condition in \eqref{STABCOND} requires that 0 is the only eigenvalue with non-negative real part, and the second condition requires that the algebraic multiplicity of 0 is 1. In other words, condition \eqref{STABCOND} requires that the linearized dynamics contracts (with exponential rate) any transverse component to the PLS circle. As such, this condition appears to be the minimal condition to ensure asymptotic nonlinear stability of the PLS circle itself.

Regarding the condition $\|\widehat{g} \|_a < +\infty$, since $g$ is real-valued, we have $\widehat{g}(-\tau) = \overline{\widehat{g}(\tau)}$. Hence $\|\widehat{g} \|_a < +\infty$ automatically yields $\int_\R e^{2a |\tau|} |\widehat{g}(\tau)|^2 d\tau < +\infty$, and then the Paley-Wiener Theorem implies that $g$ must be analytic in a horizontal strip  around the real axis.

In addition, Theorem \ref{thmNL} has a direct consequence on the behavior of the order parameter in the original Kuramoto model in finite dimension. Of note, PLS may also exist associated with discrete frequency marginals, and their expression remains the same, up to substitution of $g(\omega)d\omega$ by an arbitrary discrete probability $g(d\omega)$. However, since the Fourier transform of a discrete measure does not decay at all \cite{L70}, the first condition of Theorem \ref{thmNL} cannot hold in this case. In fact, we suspect that asymptotic stability does not hold in finite dimension.

Instead, one can use the continuous dependence of solutions of the Kuramoto equation on initial conditions \cite{L05},
to infer some control on arbitrary large time scales. In order to quantify this property, let $d_\text{BL}(\cdot,\cdot)$ be the bounded Lipschitz (or Monge-Kantorovich with exponent 1) distance on probability measures over the cylinder. Theorem 1 in \cite{L05} states the existence of $\gamma\in\R^+$ such that the solutions $t\mapsto f(t)$ and $t\mapsto f'(t)$ of the Kuramoto equation with initial condition $f_\text{in}$ and $f'_\text{in}$ respectively, satisfy
\[
d_\text{\rm BL}(f'(t),f(t))\leq d_\text{\rm BL}(f'_\text{in},f_\text{in})e^{\gamma t},\ \forall t\in\R^+.
\]
Using also the Sobolev embedding-based control of $\widehat{f}_1(t,0)=\overline{r(t)}$ by $\, \|\widehat{f(t)}\|_{a,0}$ (see details in the proof of Lemma \ref{WPKF} below), the following comment immediately results from Theorem \ref{thmNL}.
\begin{Cor}
Under the conditions of Theorem \ref{thmNL}, there exists $C_1, C_2 > 0$  such that the order parameter $r'(t)$ of the solution of the Kuramoto equation with initial data $f'_\text{\rm in}$ (of frequency marginal $g'(d\omega)$) satisfies the following inequality
\[
\left|r'(t)-r_\text{s}e^{-i\Theta_\infty} \right| \leq C_1 \, e^{-bt}+ C_2 \, d_\text{\rm BL}(f'_\text{\rm in},f_\text{\rm in})e^{\gamma t},\ \forall t\in \R^+,
\]
where $f_\text{\rm in}$ and $\Theta_\infty\in\T$ are as in Theorem \ref{thmNL}.
\end{Cor}
In particular if, in a finite dimensional system with initial empirical measure $f'_\text{in}$, the initial distance $|r'(0)-r_\text{s}e^{-i\Theta_\infty}|$ is large but the distance $d_\text{\rm BL}(f'_\text{in},f_\text{in})$ to an initial condition in the basin of attraction of the PLS circle is small (which requires that the number of oscillators be large), this statement ensures some damping of the order parameter over a large time interval.

\section{The Cauchy problem for the norm $\|\cdot\|_{a,0}$}\label{S-WELLPOSED}
The stability analysis of PLS relies on the $\|\cdot\|_{a,0}$-norm of
the (restriction of the) Fourier transform of probability
measures. However, given an arbitrary measure, its
$\|\cdot\|_{a,0}$-norm needs not to be finite. Hence, prior to the
stability analysis, we need to make sure that the set of measures with
finite $\|\cdot\|_{a,0}$-norm is invariant under the Kuramoto
flow. This is the purpose of this section, whose main result is the
following statement.
\begin{Pro} \label{WP}
Assume that $\|\widehat{g}\|_a <+\infty$ for some $a>0$ and let $t\mapsto f(t)$ be a solution of the Kuramoto equation. If  $\|\widehat{f(0)}\|_{a,0} < +\infty$, then
\[
\sup_{t\in [0,T]}\|\widehat{f(t)}\|_{a,0} <+\infty, \quad \text{for all } T > 0.
\]
Moreover, the map $t\mapsto \widehat{f(t)}$ is strongly continuous.
\end{Pro}
Before proving this statement, we express the Kuramoto dynamics in Fourier variables. Given a solution $t\mapsto f(t)$ of the Kuramoto equation, the Fourier transform $u=\widehat{f(t)}$ satisfies
\begin{equation}
\partial_t u_\ell(\tau)=\ell \partial_\tau u_\ell(\tau)+\frac{K\ell}2\left(u_1(0)u_{\ell-1}(\tau)-\overline{u_1(0)}u_{\ell+1}(\tau)\right),\ \forall (\ell,\tau)\in\N\times \R,
\label{KF}
\end{equation}
with the identification $u_0=\widehat{g}$ (and where $\overline{u_1(0)}$ is nothing but the order parameter $r$).

Given $a>0$ and $k\in \R$, consider the Hilbert space
\begin{equation} \label{Xak}
{\cal X}_{a,k}  = \{ u:\N\times\R\to \C\ \text{such that}\ \| u \|_{a,k} < +\infty\}.
\end{equation}
(We obviously have ${\cal X}_{a,k_1}\subset {\cal X}_{a,k_2}$ when $k_1>k_2$.) The proof of Proposition \ref{WP} relies on the following statement whose proof is given below.
\begin{Lem}\label{WPKF}
Assume that $\|\widehat{g}\|_a <+\infty$ for some $a > 0$. For every $u_\text{\rm in}$ in ${\cal X}_{a,0}$, there exists a unique weak solution $u$ of equation \eqref{KF} that satisfies $u(0)=u_\text{\rm in}$ and
$$u \in L^\infty(0,T, {\cal X}_{a,0}) \cap L^2(0,T, {\cal X}_{a,\frac{1}{2}})$$
for all $T > 0$. Moreover, $u$ belongs to $C(\R^+, {\cal X}_{a,0})$.
\end{Lem}
{\sl Proof of Proposition \ref{WP}.} The proposition is a simple consequence of Lemma \ref{WPKF} and Theorem 15 in \cite{D14}. This theorem provides uniqueness of the solution of  \eqref{KF} in a very large class, namely among all functions $u$ satisfying
\begin{equation} \label{criterion_uniqueness}
 \sup_{t \in [0,T]} \sup_{\ell \in \N} e^{-\beta \ell} \min(1,e^{2a \tau}) | u_\ell(t,\tau) | < +\infty
 \end{equation}
for all $T > 0$,  and for some $\beta \geq 0$.

Let $f = f(t)$ be a  solution of the Kuramoto equation, starting from an initial data $f_\text{in}$, such that  $u_\text{in} := \widehat{f_\text{in}}\vert_{\N \times \R}$ satisfies $\| u_\text{in} \|_{a,0} < +\infty$. Let  $u_f := \widehat{f}\vert_{\N \times \R}$. Clearly, $u_f$ solves \eqref{KF}, and as the Fourier transform of a measure, it is uniformly bounded in $(\ell,\tau)$.  Let now $u$ be the solution of \eqref{KF} given by Lemma \ref{WPKF}, with  initial data $u_\text{in}$.  Then, $u_f$ and $u$ both satisfy the criterion \eqref{criterion_uniqueness} for any $\beta > 0$. This implies that $u_f = u$  by Theorem 15 in \cite{D14}. \hfill $\Box$

{\sl Proof of Lemma \ref{WPKF}.}
The proof of existence proceeds via an approximation scheme and a standard compactness argument based on Aubin-Lions Lemma.  We start with the  {\sl a priori} estimates that are crucial for the limit processes. (NB: these estimates are well-defined for those $\{u_\ell(\tau)\}$ that are finite vectors of smooth functions with compact support).

The first estimate is obtained by testing \eqref{KF} against $e^{2a\tau} \frac{\overline{u_\ell(\tau)}}{\ell}$. After integration in $\tau$, summation in $\ell$, and taking the real part, we obtain
\begin{align*}
& \frac{1}{2} \frac{d}{dt}\sum_{\ell\in\N}\int_\R e^{2a\tau} \frac{|u_\ell(\tau)|^2}{\ell}d\tau
+ a \sum_{\ell\in\N}\int_\R e^{2a\tau} |u_\ell(\tau)|^2  d\tau  \\
&   = - K\text{Re}\left( \sum_{\ell\in\N} \int_\R e^{2a\tau} u_1(0) u_{\ell-1}(\tau)\overline{u_\ell(\tau)}  d \tau- \sum_{\ell\in\N} \int_\R e^{2a\tau}\overline{u_1(0)} u_{\ell+1}(\tau)\overline{u_\ell(\tau)}  d \tau\right)\\
&=-K\text{Re}\left(u_1(0)\int_\R e^{2a\tau} \widehat{g}(\tau)u_1(\tau)d\tau\right),
\end{align*}
where the last equality follows from a change $\ell\mapsto \ell+1$ of index in the first sum. Using the Cauchy-Schwarz inequality, it follows that we have
\[
\frac{1}{2} \frac{d}{dt}\sum_{\ell\in\N}\int_\R e^{2a\tau} \frac{|u_\ell(\tau)|^2}{\ell}d\tau
+ a \sum_{\ell\in\N}\int_\R e^{2a\tau}| u_\ell(\tau)|^2  d\tau\leq K|u_1(0)|\|\widehat{g}\|_a \left(\int_\R e^{2a\tau} |u_1(\tau)|^2d\tau\right)^{\frac12}.
\]
Proceeding similarly for the derivative $\partial_\tau u_\ell(\tau)$ and combining the resulting inequality with the one here then yields
\[
\frac{d}{dt}\|u\|_{a,-\frac12}^2+ 2a \|u\|_{a,0}^2\leq 2\sqrt{2} K  |u_1(0)| \|\widehat{g}\|_a \|u_1\|_a,
\]
Now, using the Sobolev embedding $H^1(0,1))\hookrightarrow C([0,1])$, we infer
\begin{equation}
|u_1(0)|\leq C\|u\|_{a,-\frac12},
\label{SOBEMB}
\end{equation}
for some $C\in\R^+$. We also have $\|u_1\|_a\leq \|u\|_{a,-\frac12}$ and the Gronwall's Lemma and the assumption $\|\widehat{g}\|_a < +\infty$ imply the existence of $C_1\in\R^+$ such that
\[
\|u(t)\|_{a,-\frac12}^2+2a\int_0^t\|u(s)\|_{a,0}^2ds\leq  e^{C_1t}\|u_\text{in}\|_{a,-\frac12}^2,\ \forall t\in\R^+.
\]
In particular, \eqref{SOBEMB} implies
\[
\sup_{t\in [0,T]}|u_1(t,0)|<+\infty,\ \forall T\in\R^+,
\]
provided that $\|u_\text{in}\|_{a,0}<+\infty$.

With this control on $|u_1(t,0)|$ provided, we can now pass to the estimate on $\|u\|_{a,0}$. To that goal, we test \eqref{KF} against $e^{a \tau}\overline{u_\ell(\tau)}$. Proceeding similarly to as before, we obtain
\begin{align*}
& \frac{1}{2} \frac{d}{dt}\sum_{\ell\in\N}\int_\R e^{2a\tau} |u_\ell(\tau)|^2d\tau
+ a \sum_{\ell\in\N}\int_\R e^{2a\tau} \ell|u_\ell(\tau)|^2  d\tau  \\
&   = - K\text{Re}\left( \sum_{\ell\in\N} \int_\R  e^{2a\tau}  \ell u_1(0) u_{\ell-1}(\tau)\overline{u_\ell(\tau)}  d \tau- \sum_{\ell\in\N} \int_\R  e^{2a\tau} \ell \overline{u_1(0)} u_{\ell+1}(\tau)\overline{u_\ell(\tau)}  d \tau\right)\\
&=- K\text{Re}\left(u_1(0)\sum_{\ell\in\N}\int_\R  e^{2a\tau}  u_{\ell-1}(\tau)\overline{u_\ell(\tau)}  d \tau\right)\\
&\leq K|u_1(0)|\left(\|\widehat{g}\|_a\left(\int_\R e^{2a\tau} |u_1(\tau)|^2d\tau\right)^{\frac12}+\sum_{\ell\in\N}\int_\R e^{2a\tau} |u_\ell(\tau)|^2d\tau\right)
\end{align*}
Repeating the argument for the derivative $\partial_\tau u_\ell(\tau)$ then yields
\[
\frac{d}{dt}\|u\|_{a,0}^2 + 2a\|u\|_{a,\frac{1}{2}}^2\leq 2K\left(\sqrt{2}|u_1(0)|\|\widehat{g}\|_{a}\|u_1\|_{a}+|u_1(0)|\|u\|_{a,0}^2\right).
\]
Finally, we use on one hand the bound \eqref{SOBEMB} with $\|u\|_{a,0}$ instead of $\|u\|_{a,-\frac12}$ and the inequality $\|u_1\|_a \leq \|u\|_{a,0}$, and on the other hand the first estimate on $\sup_{t\in [0,T]}|u_1(t,0)|$, to conclude the existence of $C_T<+\infty$ (growing at most exponentially with $T \in \R^+$) such that the following inequality holds
\begin{equation}
\|u(t)\|_{a,0}^2+2a\int_0^t\|u(s)\|_{a,\frac{1}{2}}^2ds\leq  e^{C_Tt}\|u_\text{in}\|_{a,0}^2, \ \forall t\in [0,T].
\label{estimate_lemma}
\end{equation}
This estimate allows one to construct a global  weak solution using standard arguments. For instance, one can consider a sequence of  approximate systems, by projecting equation \eqref{KF} onto a finite number of modes:
\begin{equation}
\partial_t u^n_\ell(\tau)=    \, \ell \partial_\tau u^n_\ell(\tau)+ \mathbb{P}_n \frac{K\ell}2\left(u^n_1(0)u^n_{\ell-1}(\tau)-\overline{u^n_1(0)}u^n_{\ell+1}(\tau)\right),\ \forall (\ell,\tau)\in \{1,\cdots,n\} \times \R.
\label{KFn}
\end{equation}
Here, $\mathbb{P}_n$ is the projection onto modes $\ell\in\{1,\cdots, n\}$. The
approximate initial data $u^n(0) := u^n_\text{in}$ is taken
smooth, zero for $\ell > n$ and $|\tau| > n$, and such that it converges to
$u_\text{in}$ in ${\cal X}_{a, 0}$. For any given $n$, \eqref{KFn} is a simple
transport equation with a smooth semilinear term and a smooth and
compactly supported initial data.  The existence of a local in time
solution $u^n$ is well-known \cite{R12}. The solution is smooth and compactly
supported, with $\text{supp}(u(t)) \subset \{1,\cdots,n\} \times [-n(1+t),n-t]$.  Moreover, the previous
{\it a priori} estimates extend straightforwardly to this approximate
equation, {\sl viz.}
\begin{equation*}
\|u^n(t)\|_{a,0}^2+2a\int_0^t\|u^n(s)\|_{a,\frac{1}{2}}^2ds\leq  e^{C_{T} t}\|u^n_\text{in}\|_{a,0}^2 \leq C'  e^{C_T t}, \ \forall t\in [0,T],
\end{equation*}
for any $T$ less than the maximal time of existence  $T^n$. It follows  in particular that $T^n$  is infinite. Indeed, assume {\it a contrario} that  $T^n$ is finite. As $u^n$ is compactly supported, the previous bound implies   that $u_n$ belongs to  $L^\infty((0,T^n) \times \{1,\cdots,n\} \times \R)$. This prevents blow up of the solution in finite time, and we get a contradiction.

Let $T > 0$. From the bound on $(u^n)_{n \in \N}$ in
$L^\infty(0,T,{\cal X}_{a,0})$, one can obtain a bound on the sequence
$(\pa_t u^n)_{n \in \N}$, using equation \eqref{KFn}. More precisely, the quantity
$ h^n_\ell(t,\tau) := \frac{u^n_\ell(t,\tau)}{\ell}$ is such that
\[
  (\pa_t h^n)_{n \in \N} \:  \text{ is bounded in } \:  L^\infty(0,T,
  \ell^2(\N , L^2(e^{2a \tau} d\tau)).
\]
Thus, $(u^n_1)_{n \in \N}$ is bounded in
$\displaystyle L^\infty(0,T, H^1(-1,1))$ and $(\pa_t u^n_1)_{n \in \N}$
is bounded in $\displaystyle L^\infty(0,T,L^2(-1,1))$. By Aubin-Lions Lemma, one
obtains the strong convergence of a subsequence of
$\displaystyle (u^n_1(\cdot,0))_{n \in \N}$ in
$L^\infty(0,T)$. Together with the weak compactness of $(u^n)_{n \in \N}$ in
$L^\infty(0,T , {\cal X}_{a,0}) \cap L^2(0,T , {\cal
  X}_{a,\frac{1}{2}})$,
this allows to take the limit $n \rightarrow +\infty$ in \eqref{KFn}
and yields the existence of a solution $u$ of \eqref{KF}.

For the proof of uniqueness, we use an energy estimate for the difference $v=u_2-u_1$ of two solutions. Proceeding similarly to as for the {\sl a priori} estimate above, one first obtains
\begin{align*}
& \frac{1}{2} \frac{d}{dt}\sum_{\ell\in\N}\int_\R e^{2a \tau}\frac{|v_\ell(\tau)|^2}{\ell^2}d\tau  + a \sum_{\ell\in\N}\int_\R e^{2a \tau} \frac{|v_\ell(\tau)|^2}{\ell}  d\tau =\\
&- \frac{K}2\text{Re}\left((u_1)_1(0)\sum_{\ell\in\N}\int_\R e^{2a \tau} \frac{v_{\ell}(\tau)}{\ell}\frac{\overline{u_{\ell+1}(\tau)}}{\ell+1}  d \tau + v_1(0)\sum_{\ell\in\N}\int_\R  e^{2a\tau} (u_2)_{\ell-1}(\tau)\frac{\overline{v_{\ell}(\tau)}}{\ell}  d \tau\right.\\
&\left. -\overline{v_1(0)}\sum_{\ell\in\N}\int_\R  e^{2a\tau} (u_2)_{\ell+1}(\tau)\frac{\overline{v_{\ell}(\tau)}}{\ell}  d \tau\right)
\end{align*}
and then
\[
\frac{d}{dt}\|v\|_{a,-1}^2+2a\|v\|_{a,-1/2}^2\leq C'\|v\|_{a,-1}^2,
\]
for some $C'\in\R^+$. Applying Gronwall's Lemma, the assumption $v(0)=0$ implies that $v(t)=0$ for all $t>0$ as desired.

To prove continuity, letting $h_\ell(t,\tau)=\frac{u_\ell(t,\tau)}{\ell}$, we first observe that
\[
h\in L^\infty\left(0,T,{\cal X}_{a, 1}\right)\ \text{and}\ \partial_t h\in L^\infty\left(0,T,\ell^2(\N,L^2(e^{2a\tau} d\tau))\right),\ \forall T\in\R^+.
\]
From standard functional analysis (see {\sl e.g.}\ Theorem 2.1 in \cite{S66}), it follows that $h$ is weakly continuous in time  with values in ${\cal X}_{a, 1}$, and thus $u \in C_w(\R^+, {\cal X}_{a, 0})$. Moreover, since ${\cal X}_{a, 0}$ is a Hilbert space (hence a uniformly convex space) to obtain strong continuity, it suffices to prove that
\[
t\mapsto \|u(t)\|_{a,0}
\]
is continuous (see {\sl e.g.}\ Proposition 3.32 in \cite{B11}). We consider separately the cases $t=0+$ and $t>0$.

Right continuity at 0 is rather straightforward. On one hand weak continuity implies
\[
\liminf_{t\to 0^+}\|u(t)\|_{a,0}\geq \|u(0)\|_{a,0}.
\]
On the other hand, the estimate \eqref{estimate_lemma} above implies
\[
\limsup_{t\to 0^+}\|u(t)\|_{a,0}\leq \|u(0)\|_{a,0}.
\]

For $t>0$, we use the regularization effect induced by the
weight. The integral term of \eqref{estimate_lemma} shows that
\[
\|u(\delta)\|_{a,\frac{1}{2}}<+\infty,\ \text{for a.e.}\ \delta \in\R^+.
\]
Take any such $\delta$. By mimicking the arguments above, one can construct a solution $\tilde u = \{\tilde u_\ell(t,\tau)\}$ of \eqref{KF} over $(\delta,+\infty)$ satisfying
$$ \tilde u \in L^\infty(\delta,T , {\cal X}_{a, \frac{1}{2}}) \cap L^2(\delta,T,  {\cal X}_{a,1}), \quad \text{for all } T > \delta,  $$
with $\tilde u(\delta) = u(\delta)$. By invoking the uniqueness of the solution in $L^\infty(\delta,T, {\cal X}_{a, 0}) \cap   L^2(\delta,T, {\cal X}_{a, \frac{1}{2}})$, we deduce $\tilde u = u$. It follows in particular that $u  \in L^2(\delta, T,{\cal X}_{a, 1})$ for any $T > \delta > 0$. It is then easily seen that \eqref{KF} reads
$$ \pa_t u_l(\tau) - l \pa_\tau u_l(\tau)  = F_l(\tau) $$
where $F \in L^2(\delta, T,{\cal X}_{a, 0})$ for any $T > \delta > 0$. Using the explicit formula
$$ u_\ell(t,\tau) = u_\ell(\delta,\tau) +  \int_\delta^t F_\ell(s,\tau + \ell(t-s)) ds, \quad t > \delta > 0, $$
one can check that $u$ is continuous at positive times with values in ${\cal X}_{a, 0}$. \hfill $\Box$

\section{Spectral analysis of the linearized dynamics}\label{S-LINEAR}
As for the Cauchy problem, the analysis of the perturbation dynamics proceeds in Fourier variables. Assuming a PLS $f_\text{s}$ and inserting the expression $\widehat{f_\text{s}}+u$ in the Fourier formulation \eqref{KF} of the Kuramoto equation, the time evolution of the perturbation $u=\{u_\ell(\tau)\}_{\N\times\R}$ (NB: from now on, $u$ denotes a perturbation to the PLS $\widehat{f_\text{s}}$) turns out to be governed by the equation
\begin{equation} \label{KF_condensed}
\partial_t u=Lu+Qu,
\end{equation}
where $L=L_1+L_2$ and for all $(\ell,\tau)\in\N\times\R$, we have, using the notation $u_0(\tau)=0$,
\[
(L_1u)_\ell(\tau)=\ell\left(\partial_\tau u_\ell(\tau)+\frac{K r_\text{s}}2\left(u_{\ell-1}(\tau)-u_{\ell+1}(\tau)\right)\right),
\]
and
\[
(L_2u)_\ell(\tau)=\frac{K \ell}2\left(u_1(0)(\widehat{f_\text{s}})_{\ell-1}(\tau)-\overline{u_1(0)}(\widehat{f_\text{s}})_{\ell+1}(\tau)\right),
\]
and the operator $Q$ collects the nonlinear terms
\[
(Qu)_\ell(\tau)=\frac{K \ell}2\left(u_1(0)u_{\ell-1}(\tau)-\overline{u_1(0)}u_{\ell+1}(\tau)\right).
\]

This section deals with the analysis of linear terms. Of note, while the operator $L$ is $\R$-linear, it is not $\C$-linear simply because $L_2$ does not satisfy this property. In order to get a $\C$-linear operator and to investigate its spectral properties, one may consider the real and imaginary parts separately, as in \cite{MS07,OW13}. We use an alternative approach here, based on complex conjugates. Given $u=\{u_\ell(\tau)\}_{\N\times\R}$ and $v=\{v_\ell(\tau)\}_{\N\times\R}$ (which is a substitute for $\bar{u}$), let
\[
\text{\sl u}=\{\text{\sl u}_\ell(\tau)\}_{\N\times\R}\ \text{where}\ \text{\sl u}_\ell(\tau)={u_\ell(\tau) \choose v_\ell(\tau)}\in\C^2,\ \forall (\ell,\tau)\in\N\times\R,
\]
and consider the operator ${\cal L}={\cal L}_1+{\cal L}_2$ defined by
\[
({\cal L}_1\text{\sl u})_\ell(\tau)={(L_1u)_\ell(\tau)\choose (L_1v)_\ell(\tau)}
\]
and
\[
({\cal L}_2\text{\sl u})_\ell(\tau)=\frac{K \ell}2
\left(\begin{array}{cc}(\widehat{f_\text{s}})_{\ell-1}(\tau)&-(\widehat{f_\text{s}})_{\ell+1}(\tau)\\
-\overline{(\widehat{f_\text{s}})_{\ell+1}(\tau)}&\overline{(\widehat{f_\text{s}})_{\ell-1}(\tau)}\end{array}\right){u_1(0)\choose v_1(0)}.
\]
The operators ${\cal L}_i$ are defined in such a way that when $v_\ell(\tau)=\overline{u_\ell(\tau)}$, we have
\[
({\cal L}_i\text{\sl u})_\ell(\tau)={(L_iu)_\ell(\tau)\choose \overline{(L_iu)_\ell(\tau)}},\ \text{for}\ i=1,2.
\]
Given $a>0$ and $k\in\R$, let
\[
D_{a,k}=\left\{u\in {\cal X}_{a,k}\ :\ L_1u\in {\cal X}_{a,k}\right\},
\]
be the domain of the operator $L_1$. Thanks to Proposition \ref{UNIQPLS} in Appendix \ref{A-STAB}, $L_2$ is also well-defined on $D_{a,k}$; hence the product $D_{a,k}^2$ is a domain for the operator ${\cal L}$ on ${\cal X}_{a,k}$. Since they contain finite vectors of smooth functions with compact support, the domains $D_{a,k}$ are dense in their respective space.

The results of this section are collected in the following statement.
\begin{Pro}
Assume that $\|\widehat{g}\|_{a}<+\infty$ for some $a>0$ and assume that a PLS $f_\text{s}$ with marginal density $g$ and order parameter $r_\text{s}\in\R^+$ exists. The corresponding operator ${\cal L}$ has the following properties on ${\cal X}_{a,k}^2$ for $k\in\{-1,0\}$.
\begin{itemize}
\item[$\bullet$] it generates a $C^0$-semigroup,
\item[$\bullet$] its essential spectrum lies in the half-plane $\text{\rm Re}(\lambda)\leq -a$,
\item[$\bullet$] for every $\epsilon>0$, its spectrum in the half-plane $\text{\rm Re}(\lambda)>-a+\epsilon$ consists of finitely many eigenvalues with finite multiplicity; 0 is always one of them,
\item[$\bullet$] let $x>-a$ be such that the line $\{\lambda \in \C: \text{\rm Re}(\lambda)= x\}$  does not contain any eigenvalue. Then we have
 \begin{equation*}
    \sup_{y\in\R}
    \| ((x+iy) \text{\rm Id} - {\cal L})^{-1} \|_{{\cal X}_{a,-1} \to
      {\cal X}_{a,0}} < +\infty.
  \end{equation*}
\end{itemize}
\label{thmL}
\end{Pro}
The proof is given in the two sections below. First, we obtain estimates for the resolvent of $L_1$ (which is $\C$-linear) from where strong continuity of the semi-group readily follows. Then, the spectrum is described and the last estimate is shown.

Proposition \ref{thmL} suggests that a condition for asymptotic stability of the circle $\{R_\Theta f_\text{s}\}$ of PLS is that the spectrum of ${\cal L}$ in the half-plane $\text{Re}(\lambda)\geq 0$ consists of the sole eigenvalue 0, and this eigenvalue is simple. (The analysis of nonlinear terms in section \ref{S-NLSTAB} below shows that this is indeed the case.)
In principle, it would suffice to impose stability of the restriction of ${\cal L}$ to the subspace where $v_\ell(\tau)=\overline{u_\ell(\tau)}$. However, this consideration does not make any difference. We shall verify that the complexification process does not introduce any unstable spurious mode, namely that ${\cal L}$ does not have unstable spectrum when the original linear equation $\partial_t u=Lu$ is stable (see end of the section on the spectral analysis).

\subsection{Resolvent estimates for $L_1$}
\begin{Lem}
\label{thm:resolvent-l1-generation}
Let $k\in \{-1,0\}$. The resolvent set of $L_1$ over ${\cal X}_{a,k}$ contains the half-plane $\text{\rm Re}(\lambda)>-a$ and we have
$$
\| (\lambda \text{\rm Id}-L_1)^{-1}\|_{{\cal X}_{a,-1} \to {\cal X}_{a,0}}\leq \frac1{\min\{a,\text{\rm Re}(\lambda)+a\}},\ \forall \lambda\in\C\ :\ \text{\rm Re}(\lambda)>-a.
$$
Moreover, letting $\lambda_0=-a+\frac{Kr_\text{\rm s}}2$, we have in the half-plane $\text{\rm Re}(\lambda)>\lambda_0$
\[
\|(\lambda \text{\rm Id}-L_1)^{-1}\|_{a,k}\leq \frac1{\text{\rm Re}(\lambda)-\lambda_0}.
\]
\end{Lem}
Clearly, $L_1$ is a closed operator on $D_{a,k}$; hence applying the Hille-Yosida Theorem with the second estimate in the Lemma implies that $L_1$ generates a $C^0$-semigroup. The same property certainly holds for the operator ${\cal L}_1$. Moreover, thanks to the property $\|\widehat{f_\text{s}}\|_{a,k+\tfrac12}<+\infty$ (see Proposition \ref{UNIQPLS} in Appendix \ref{A-STAB}), the perturbation ${\cal L}_2$ is bounded on ${\cal X}_{a,k}^2$; hence, ${\cal L}$ also generates a $C^0$ semigroup (see {\sl e.g.}\ \cite{K95}). The first item of Proposition \ref{thmL} is proved.

In addition, the first estimate of Lemma \ref{thm:resolvent-l1-generation} implies that  $\| (\lambda \text{\rm Id}-L_1)^{-1}\|_{a,k}$ is uniformly bounded over  the half-plane $\text{Re}(\lambda)> -a + \eps$ for any $\eps > 0$.  The Gearhart-Pr\"uss Theorem (see {\sl e.g.}\ Corollary 2.2.5 in \cite{vN96}) then implies that the semigroup $e^{tL_1}$ must be exponentially stable, more precisely that there exist $b>a$ and $C \in\R^+$ such that
\begin{equation}
\|e^{tL_1}\|_{a,k}\leq C e^{-bt},\ \forall t\in\R^+.
\label{STABESTIMATE}
\end{equation}

{\sl Proof of Lemma \ref{thm:resolvent-l1-generation}.} To derive the claimed inequalities,  we consider the resolvent equation
\[
(\lambda \text{\rm Id}-L_1)u=v.
\]
The second estimate is obtained by testing against $e^{2a\tau}\ell^{2k}\overline{u_\ell(\tau)}$, under the assumption $\|u\|_{a,k}<+\infty$. After integration in $\tau$, summation in $\ell$ and taking the real part, we obtain
\begin{align*}
&\text{Re}(\lambda)\sum_{\ell\in\N}\int_\R e^{2a\tau}\ell^{2k}|u_\ell(\tau)|^2d\tau +a \sum_{\ell\in\N}\int_\R e^{2a\tau}\ell^{2k+1}|u_\ell(\tau)|^2d\tau \\
&-\frac{K r_\text{s}}2\text{Re}\left(\sum_{\ell\in\N}\int_\R e^{2a\tau}\ell^{2k+1}(u_{\ell-1}(\tau)\overline{u_\ell(\tau)}-\overline{u_{\ell+1}(\tau)}u_{\ell}(\tau))d\tau\right)\leq \|u\|_{a,k}\|v\|_{a,k}
\end{align*}
A change $\ell\mapsto \ell+1$ of index in the third sum yields, also using $u_0(\tau)=0$ and simplifying the expression of $(\ell+1)^{2k+1}-\ell^{2k+1}$ for $k\in \{-1,0\}$
\begin{align*}
&\left(\text{Re}(\lambda)+a\right)\sum_{\ell\in\N}\int_\R e^{2a\tau}\ell^{2k}|u_\ell(\tau)|^2d\tau\\
&\leq \frac{K r_\text{s}}2\text{Re}\left(\sum_{\ell\in\N}\int_\R e^{2a\tau}\left((\ell+1)^{2k+1}-\ell^{2k+1}\right)u_{\ell}(\tau)\overline{u_{\ell+1}(\tau)}d\tau\right)+ \|u\|_{a,k}\|v\|_{a,k}\\
&\leq \frac{K r_\text{s}}2\|u\|_{a,k}^2+ \|u\|_{a,k}\|v\|_{a,k}.
\end{align*}
As in the proof of Lemma \ref{WPKF}, one can proceed similarly for the derivative $\partial_\tau u_\ell(\tau)$ and the second estimate easily follows.

For the first estimate, we proceed similarly for $k=-\frac12$ and use the inequality
\[
\text{Re}(\lambda)\sum_{\ell\in\N}\int_\R e^{2a\tau}|u_\ell(\tau)|^2d\tau\leq \text{Re}(\lambda)\sum_{\ell\in\N}\int_\R e^{2a\tau}\ell^{-1}|u_\ell(\tau)|^2d\tau,\ \forall \lambda\ :\ \text{Re}(\lambda)\leq 0,
\]
to obtain
\[
\min\{a,\text{\rm Re}(\lambda)+a\}\sum_{\ell\in\N}\int_\R e^{2a\tau}|u_\ell(\tau)|^2d\tau\leq \|u\|_{a,0}\|v\|_{a,-1},
\]
and then
\[
\|u\|_{a,0}\leq \frac{\|v\|_{a,-1}}{\min\{a,\text{\rm Re}(\lambda)+a\}},\ \forall \lambda\in\C\ :\ \text{\rm Re}(\lambda)>-a,
\]
from where the first estimate follows suit.
Finally, standard arguments ({\sl e.g.}\ Galerkin approximation) based on this estimate show that $\lambda \text{\rm Id}-L_1$ is invertible for $\text{\rm Re}(\lambda)>-a$. \hfill $\Box$

\subsection{The spectrum of ${\cal L}$}
Lemma \ref{thm:resolvent-l1-generation} obviously implies that the spectrum of ${\cal L}_1$ must be contained in the half-plane $\text{Re}(\lambda)\leq -a$. Moreover, the perturbation ${\cal L}_2$ has finite rank; hence the essential spectrum of ${\cal L}$ must be contained in the same region, as claimed in the second item of Proposition \ref{thmL}.

To prove the third item of that Proposition, we characterize the eigenvalues in the complementary region $\text{Re}(\lambda)> -a$. To that goal, consider the $2\times 2$ matrix defined by
\[
M'(\lambda,r_\text{s})=\left(\begin{array}{cc}
((\lambda \text{Id}-L_1)^{-1}\widehat{p_0})_1(0)&-((\lambda \text{Id}-L_1)^{-1}\widehat{p_2})_1(0)\\
-\overline{((\overline{\lambda} \text{Id}-L_1)^{-1}\widehat{p_2})_1(0)}&\overline{((\overline{\lambda} \text{Id}-L_1)^{-1}\widehat{p_0})_1(0)}\end{array}\right),
\]
where the vector $\widehat{p_k}$ ($k\in\Z$) is defined by
\[
(\widehat{p_k})_\ell(\tau)=\ell(\widehat{f_\text{s}})_{\ell+k-1}(\tau),\ \forall (\ell,\tau)\in \N\times\R.
\]
\begin{Lem}
The number $\lambda$ with $\text{Re}(\lambda)> -a$ is an eigenvalue of ${\cal L}$ iff $\text{\rm det}\left(\text{\rm Id}-\frac{K}2M'(\lambda,r_\text{s})\right)=0$.
\label{QUALIFEIGEN}
\end{Lem}
In addition, that there are finitely many eigenvalues, with finite
multiplicity, in every half-plane
$\text{Re}(\lambda)>-a+\epsilon$ ($\epsilon>0$) is a
consequence of the holomorphic dependence on $\lambda$ and the following limits
\[
\lim_{x\to +\infty}\sup_{y\in\R}(((x+iy) \text{Id}-L_1)^{-1}\widehat{p_k})_1(0)=0\ \text{and}\
\lim_{y\to\pm \infty}\sup_{x\geq -a+\eps}(((x+iy) \text{Id}-L_1)^{-1}\widehat{p_k})_1(0)=0.
\]
The first limit is a direct consequence of Lemma \ref{thm:resolvent-l1-generation}. For the second limit, we use the expression of the resolvent as the Laplace transform of the semi-group (see e.g.\ Section 1.7 in \cite{P83}) to obtain
\[
(((x+iy) \text{Id}-L_1)^{-1}\widehat{p_k})_1(0)=\int_{\R^+}e^{-iyt}(e^{-x t}e^{tL_1}\widehat{p_k})_1(0)dt.
\]
Moreover, the inequality \eqref{STABESTIMATE} implies that
$t\mapsto (e^{-x t}e^{tL_1}\widehat{p_k})_1(0)$ is uniformly
absolutely integrable for $x\geq -a+\eps$. The limit then
follows from the Riemann-Lebesgue Lemma.

Finally, we observe that the rotation symmetry $R_\Theta$ of the Kuramoto equation expresses as a phase symmetry in Fourier variables, {\sl i.e.}\ if $t\mapsto u(t)=\{u_\ell(t,\tau)\}_{\N\times\R}$ satisfies \eqref{KF}, then for every $\Theta\in\T$, the trajectory $t\mapsto \widehat{R}_\Theta u(t)$, where $(\widehat{R}_\Theta u(t))_{\ell,\tau}=e^{i\Theta\ell}u_\ell(t,\tau)$, also solves that equation. As noted in \cite{MS07}, this indifference to phase changes implies that we must have
\begin{equation}
Lu=0\ \text{for}\ u=D\widehat{R}\widehat{f_\text{s}}.
\label{0EIGEN}
\end{equation}
where $D\widehat{R}:=\frac{d\widehat{R}_\Theta}{d\Theta}|_{\Theta=0}$ is the symmetry infinitesimal generator and writes
\[
(D\widehat{R}u)_\ell(\tau)=i\ell u_\ell(\tau).
\]
In particular, 0 then must be an eigenvalue of ${\cal L}$. The proof of the third item of Proposition \ref{thmL} is complete.

{\sl Proof of Lemma \ref{QUALIFEIGEN}.} Let
\[
\text{\rm U}=\{\text{\rm U}_\ell(\tau)\}_{\N\times\R}\ \text{where}\ \text{\rm U}_\ell(\tau)={U_\ell(\tau) \choose V_\ell(\tau)}\in\C^2,
\]
then the resolvent equation $(\lambda \text{Id}-{\cal L})\text{\sl u}=\text{\rm U}$ can be written in the region $\text{Re}(\lambda)> -a$
\[
\left(\text{Id}-(\lambda \text{Id}-{\cal L}_1)^{-1}{\cal L}_2\right)\text{\sl u}=(\lambda \text{Id}-{\cal L}_1)^{-1}\text{\rm U}.
\]
An important property is that the vector ${\cal L}_2\text{\sl u}$, and hence $(\lambda \text{Id}-{\cal L}_1)^{-1}{\cal L}_2\text{\sl u}$, only involves the component $\text{\sl u}_1(0)$ of $\text{\sl u}$. Using also the commutation $L_1\overline{u}=\overline{L_1u}$, it follows that the component $(\ell,\tau)=(1,0)$ of the resolvent equation writes
\begin{equation}
\left(\text{Id}-\frac{K}2M'(\lambda,r_\text{s})\right)\text{\sl u}_1(0)=(\lambda \text{Id}-{\cal L}_1)^{-1}\text{\rm U}_1(0),
\label{EQ10}
\end{equation}
Therefore, in the case where $\text{Id}-\frac{K}2M'(\lambda,r_\text{s})$ is invertible, let $\text{\sl u}^\ast=\{\text{\sl u}^\ast_\ell(\tau)\}_{\N\times\R}$ be any vector for which $\text{\sl u}^\ast_1(0)$ solves \eqref{EQ10}.  We infer that the resolvent equation has a solution given by
\[
\text{\sl u}=(\lambda \text{Id}-{\cal L}_1)^{-1}\left({\cal L}_2\text{\sl u}^\ast+\text{\rm U}\right),
\]
which is unique since $(\lambda \text{Id}-{\cal L}_1)^{-1}{\cal L}_2\text{\sl u}^\ast$ only involves the component $\text{\sl u}^\ast_1(0)$.

On the other hand, if $\text{det}\left(\text{Id}-\frac{K}2M'(\lambda,r_\text{s})\right)=0$, let $\text{\sl u}^\dag$ be with component $\text{\sl u}^\dag_1(0)\in \text{Ker}\left(\text{Id}-\frac{K}2M'(\lambda,r_\text{s})\right)$. Using once again that $(\lambda \text{Id}-{\cal L}_1)^{-1}{\cal L}_2\text{\sl u}^\dag$ only involves $\text{\sl u}^\dag_1(0)$, one directly checks that $(\lambda \text{Id}-{\cal L}_1)^{-1}{\cal L}_2\text{\sl u}^\dag$ is an eigenvector of ${\cal L}$ with eigenvalue $\lambda$.
Consequently, $\det \left(\text{\rm Id}-\frac{K}2 M'(\lambda,r_\text{s})\right)= 0$ iff $\lambda$ is an eigenvalue of ${\cal L}$ in the half-plane $\text{Re}(\lambda)>-a$. \hfill $\Box$

Finally, the last item of Proposition \ref{thmL} can be shown by combining the first claim in Lemma \ref{thm:resolvent-l1-generation} with  the expression of the resolvent in the proof of Lemma \ref{QUALIFEIGEN}.

{\bf Absence of spurious modes.} To complete this section, we show that the complexification process does not generate unstable spurious modes, {\sl i.e.}\ to any eigenvalue $\lambda$ with $\text{Re}(\lambda)>0$ of ${\cal L}$ (resp.\ non-zero eigenvalue on the imaginary axis), corresponds a diverging (resp.\ rotating) solution of $\partial_t u=Lu$. We consider the cases $\text{Im}(\lambda)\neq 0$ and $\lambda\in\R$ separately.

{\sl Case $\text{\rm Im}(\lambda)\neq 0$.} Given that $\lambda$ and $\overline{\lambda}$ both are eigenvalues of ${\cal L}$, the trajectory $t\mapsto \text{\rm U}(t)$, (uniquely) defined by
\[
\text{\rm U}(t)=e^{\lambda t}(\lambda \text{Id}-{\cal L}_1)^{-1}{\cal L}_2\text{\sl u}+e^{\overline{\lambda} t}(\overline{\lambda} \text{Id}-{\cal L}_1)^{-1}{\cal L}_2\overline{\text{\sl u}},
\]
where $\text{\sl u}$ is with component $\text{\sl u}_1(0)={u_1(0) \choose v_1(0)}\in \text{Ker}\left(\text{Id}-\frac{K}2M'(\lambda,r_\text{s})\right)$ and $\overline{\text{\sl u}}$ is with component $\overline{\text{\sl u}}_1(0)={\overline{v}_1(0) \choose \overline{u}_1(0)}\in \text{Ker}\left(\text{Id}-\frac{K}2M'(\overline{\lambda},r_\text{s})\right)$, is a solution of the equation $\partial_t \text{\sl u}={\cal L}\text{\sl u}$. Moreover, this solution components satisfy $(V(t))_\ell(\tau)=\overline{(U(t))_\ell(\tau)}$; hence the definition of ${\cal L}$ implies that $\{(U(t))_\ell(\tau)\}_{\N\times\R}$ satisfies the equation $\partial_t u=Lu$.

{\sl Case $\lambda\in\R$.} In this case, the matrix $\text{Id}-\frac{K}2M'(\lambda,r_\text{s})$ must be Hermitian and of the form
\[
\rho\left(\begin{array}{cc}
e^{i\varphi_0}&e^{i\varphi_2}\\
e^{-i\varphi_2}&e^{-i\varphi_0}
\end{array}\right),
\]
for some $\varphi_0,\varphi_2\in\T$ and $\rho\in\R^+$. Clearly, as $\rho \neq 0$, the kernel of this matrix is of the form $\text{Span}\left\{ {e^{i\varphi}\choose e^{-i\varphi}}\right\}$ for a given $\varphi\in\C$. Letting $\text{\sl u}_1(0)$ be in this kernel and $\text{\rm U}(t)=e^{\lambda t}(\lambda \text{Id}-{\cal L}_1)^{-1}{\cal L}_2\text{\sl u}$, we have that the first component $\{(U(t))_\ell(\tau)\}_{\N\times\R}$ must also satisfy the equation $\partial_t u=Lu$ in this case.

\subsection{Derivation of the stability condition}\label{S-DERIVSTAB}
As said after Proposition \ref{thmL}, the condition for linear stability of the PLS circle is that ${\cal L}$ has no eigenvalue for $\text{Re}(\lambda)\geq 0$, except 0 and this eigenvalue is simple.
Given Lemma \ref{QUALIFEIGEN}, these requirements are equivalent to the following ones
\begin{equation*}
\left\{
\begin{aligned}
& \det\left( \text{Id}-\frac{K}{2}M'(\lambda,r_\text{s}) \right) \neq 0, \ \forall \lambda \neq 0, \: \text{Re}(\lambda) \geq  0, \\
& \text{$\lambda=0$ is a simple zero of the holomorphic function $\det\left( \text{Id}-\frac{K}{2}M'(\lambda,r_\text{s}) \right)$}.
\end{aligned}
\right.
\end{equation*}
This last condition is obviously the same as  $\liminf_{\lambda \rightarrow 0} \left| \frac{1}{\lambda} \det\left( \text{Id}-\frac{K}{2}M'(\lambda,r_\text{s}) \right) \right| > 0$.

Now, condition \eqref{STABCOND} of Theorem \ref{thmNL} readily follows from the fact
that in the half-plane
$\text{Re}(\lambda)\geq 0$, the matrix $M'(\cdot,r_\text{s})$ turns
out to coincide with the matrix obtained by flipping the off-diagonal
terms in $M(\cdot,r_\text{s})$. This fact is an immediate consequence of the following statement (which also ensures that $M(\lambda,r_\text{s})$ is well-defined).
\begin{Lem}
The quantities $J_k(\lambda,r_\text{s})$ introduced in \eqref{DEFJK} are well-defined for all $k\in\N\cup\{0\}$ and $\text{\rm Re}(\lambda)\geq 0$, and we have
\[
J_k(\lambda,r_\text{s})=((\lambda \text{\rm Id}-L_1)^{-1}\widehat{p_k})_1(0).
\]
\label{EQUIVM}
\end{Lem}

{\sl Proof.} We use inverse Fourier transforms with respect to $\tau$. Let
\[
(p_k)_\ell(\omega)=\frac1{2\pi}\int_\R e^{i\omega\tau}(\widehat{p_k})_\ell(\tau)d\tau,
\]
and let $\check{L}_1$ be the inverse Fourier transform of $L_1$, {\sl i.e.}\ $\widehat{\check{L}_1p_k}=L_1\widehat{p_k}$,
when passing to Fourier transforms with respect to $\tau$. We have
\[
((\lambda \text{\rm Id}-L_1)^{-1}\widehat{p_k})_1(0)=\int_\R((\lambda \text{\rm Id}-\check{L}_1)^{-1}p_k)_1(\omega)d\omega,
\]
provided that $((\lambda \text{\rm Id}-\check{L}_1)^{-1}p_k)_1\in L^1(\R)$, and expression \eqref{FOURCOEF} in Appendix \ref{A-STAB} implies
\[
(p_k)_\ell(\omega)=\ell (\widetilde{f_\text{s}})_{\ell+k-1}(\omega)=\ell\beta^{\ell+k-1}\left(\frac{\omega}{Kr_\text{s}}\right) g(\omega).
\]
Now, using the expression
\[
(\check{L}_1u)_\ell(\omega)=\ell\left(i\omega u_\ell(\omega)+\frac{Kr_\text{s}}2(u_{\ell-1}(\omega)-u_{\ell+1}(\omega))\right),
\]
and, twice in a row, the equation \eqref{EQBETA} in Section \ref{S-STAB} that defines $\beta$, one obtains
 \[
(\check{L}_1p_k)_\ell(\omega)=-\frac{Kr_\text{s}}2\frac{1+\beta^2\left(\frac{\omega}{Kr_\text{s}}\right)}{\beta\left(\frac{\omega}{Kr_\text{s}}\right)}(p_k)_\ell(\omega)=-\left(i\omega+Kr_\text{s}\beta\left(\frac{\omega}{Kr_\text{s}}\right)\right)(p_k)_\ell(\omega),
\]
from where it results that
\[
((\lambda \text{\rm Id}-\check{L}_1)^{-1}p_k)_1(\omega)=\frac{\beta^k\left(\frac{\omega}{Kr_\text{s}}\right)g(\omega)}{\lambda +i\omega+Kr_\text{s}\beta\left(\frac{\omega}{Kr_\text{s}}\right)}.
\]
Using the expression of $\beta$ and $|\beta(\cdot)|\leq 1$, the following inequality holds
\[
\left|\frac{\beta^k\left(\frac{\omega}{Kr_\text{s}}\right)}{\lambda +i\omega+Kr_\text{s}\beta\left(\frac{\omega}{Kr_\text{s}}\right)}\right|\leq \frac{1}{\text{Re}(\lambda)},
\]
provided that $\text{Re}(\lambda)>0$. The lemma then easily follows in this case.

For $\text{Re}(\lambda)=0$, the integral defining  $J_k(\lambda,r_\text{s})$ has to be understood in a weak sense: it is  the limit as $\eps \rightarrow 0^+$ of  $J_k(\lambda + \eps,r_\text{s})$.  This limit exists because of the  continuous dependence of the resolvent $(\lambda \text{\rm Id}-L_1)^{-1}$ on $\lambda$. In practice, the value of $J_k(\lambda,r_\text{s})$ for $\lambda$ on the imaginary axis can be computed as a principal value with correction terms, using Plemelj formula as in \cite{V10}. \hfill $\Box$

\section{Proof of nonlinear stability}\label{S-NLSTAB}
With full understanding of the linearized dynamics, we can now address nonlinear terms, assuming that the stability condition \eqref{STABCOND} holds. Following a standard procedure in the stability analysis of relative equilibria (see {\sl e.g.}\ \cite{HI11}), the first step is to take advantage of the rotation symmetry $\widehat{R}_\Theta$ to get rid of the corresponding indifferent angular coordinate. Then, we shall prove asymptotic stability of the stationary state associated with the remaining variable.


\subsection{Polar-type coordinates}
To get rid of the angular coordinate, which is tangent to the PLS circle $\{\widehat{R}_\Theta \widehat{f_\text{s}}\}_{\Theta\in\T}$, we need to introduce a projection operator $P_0$ defined on ${\cal X}_{a,0}$ such that $LP_0=0$ and $P_0L=0$.

The former equality obviously implies that $P_0$ must project on $\text{Ker}(L)$. To comply with the latter equality, it is convenient to consider the analogous projection associated with the operator ${\cal L}$ on ${\cal X}_{a,0}^2$.
Let $\langle\cdot ,\cdot\rangle_{a,0}$ be the scalar product on ${\cal X}_{a,0}$ that induces the norm $\|\cdot\|_{a,0}$ and let the scalar product $\langle\cdot ,\cdot\rangle$ on ${\cal X}_{a,0}^2$ be defined by
\[
\langle\text{\sl u} ,\text{\sl u}'\rangle = \langle u,u'\rangle_{a,0}+\langle v ,v'\rangle_{a,0},\ \text{where}\ \text{\sl u}={u\choose v}\ \text{and}\ \text{\sl u}'={u'\choose v'}.
\]
The stability condition \eqref{STABCOND} implies that $\text{Ker}({\cal L})=\text{Span}{D\widehat{R}\widehat{f_\text{s}}\choose \overline{D\widehat{R}\widehat{f_\text{s}}}}$ is one-dimensional. Accordingly, let $\text{\sl u}^\ast$ be such that $\text{Ker}({\cal L}^\ast)=\text{Span}(\text{\sl u}^\ast)$ for the adjoint operator ${\cal L}^\ast$, and such that $\langle {D\widehat{R}\widehat{f_\text{s}}\choose \overline{D\widehat{R}\widehat{f_\text{s}}}},\text{\sl u}^\ast\rangle=1$. Let then the projection operator ${\cal P}_0$ be defined by
\[
{\cal P}_0\text{\sl u}=\langle  \text{\sl u},\text{\sl u}^\ast\rangle {D\widehat{R}\widehat{f_\text{s}}\choose \overline{D\widehat{R}\widehat{f_\text{s}}}},\ \forall \text{\sl u}\in {\cal X}_{a,0}^2.
\]
This operator obviously satisfies ${\cal P}_0{\cal L}\text{\sl u}=0$ for all $\text{\sl u}$. Moreover, we are going to show that $\text{\sl u}^\ast$ takes the form ${u^\ast\choose \overline{u^\ast}}$ for some $u^\ast\in {\cal X}_{a,0}$ (which then must satisfy the normalization condition $2\text{Re}\langle D\widehat{R}\widehat{f_\text{s}},u^\ast\rangle_{a,0}=1$. Therefore, we have
\[
{\cal P}_0\text{\sl u}={P_0 u\choose \overline{P_0 u}}\ \text{when}\ \text{\sl u}={u\choose\overline{u}},\
\text{where}\ P_0u=2\text{Re}\langle u,u^\ast\rangle_{a,0} D\widehat{R}\widehat{f_\text{s}}.
\]
The operator $P_0$ is the desired projection, {\sl i.e.}\ we have $LP_0u=P_0Lu=0$ for all $u\in {\cal X}_{a,0}$.

To show the promised property, let ${u\choose v}\in \text{Ker}(\cal L^\ast)$ be arbitrary. Then ${u+\overline{v}\choose v+\overline{u}}={u+\overline{v}\choose \overline{u+\overline{v}}}\in \text{Ker}(\cal L^\ast)$, which implies $\text{\sl u}^\ast={u^\ast\choose \overline{u^\ast}}$ as claimed. Indeed, either $u+\overline{v}\neq 0$ and that $\text{Ker}(\cal L^\ast)$ is one-dimensional implies that we must have  $u^\ast=\lambda (u+\overline{v})$ for some $\lambda\in\C$. Or $u+\overline{v}=0$ and then $u^\ast=\lambda i u$ for some $\lambda\in\C$.

For future purposes, we now show that
\begin{equation} \label{regularityu*}
u^* \in {\cal X}_{a,k}, \ \forall k \in \N\cup\{0\},
\end{equation}
which, in particular, implies that $P_0$ is a projection operator on every ${\cal X}_{a,k}$.

To that goal, observe that ${\cal L}^\ast$ can be explicitly computed as ${\cal L}^\ast={\cal L}_1^\ast+{\cal L}_2^\ast$ where
\begin{align*}
&{\cal L}_1^\ast\text{\sl u} = {L_1^\ast u \choose L_1^\ast v}, \ \text{where}\ (L_1^\ast u)_\ell = - \ell \left(\pa_\tau  u_\ell+2a u_\ell\right)+ \frac{K r_s}{2} \left( (\ell+1) u_{\ell+1}  - (\ell-1) u_{\ell-1} \right)
\\
&\text{and}\ ({\cal L}_2^* \text{\sl u})_\ell =  \frac{K}{2}{m_u w\delta_{\ell,1}\choose m_v w\delta_{\ell,1}}\ \text{where}\ {m_u\choose m_v}=\sum_{\ell\in\N}  \int_\R e^{ 2a \tau}\ell
\left(\begin{array}{cc}(\overline{\widehat{f_\text{s}})_{\ell-1}(\tau)} & -(\widehat{f_\text{s}})_{\ell+1}(\tau)\\
-\overline{(\widehat{f_\text{s}})_{\ell+1}(\tau)}& (\widehat{f_\text{s}})_{\ell-1}(\tau)\end{array}\right) {u_\ell(\tau) \choose v_\ell(\tau)}d\tau ,
\end{align*}
where we have used the Kronecker symbol and $w:\R\to\C$ is the function such that $\|w\|_a<+\infty$ and
\[
\langle w,w'\rangle_a=\overline{w'(0)},\ \forall w':\R\to\C\ :\ \|w'\|_a<+\infty,
\]
where $\langle \cdot,\cdot\rangle_a$ is the scalar product that generates $\|\cdot\|_a$. The existence and uniqueness of $w$ is guaranteed by the Riesz Representation Theorem, using Sobolev embedding to ensure that $w'\mapsto \overline{w'(0)}$ is a continuous linear functional.

Consequently, $u^\ast$ must satisfy the equation $(L_1^\ast u^\ast)_\ell=-c^\ast w\delta_{\ell,1}$ where
\[
c^\ast=-\sum_{\ell\in\N}  \int_\R e^{ 2a \tau}\ell\left(\overline{(\widehat{f_\text{s}})_{\ell-1}(\tau)} u_\ell^\ast(\tau) -(\widehat{f_\text{s}})_{\ell+1}(\tau)\overline{u_\ell^\ast(\tau)}\right)d\tau.
\]
Letting $x^\ast_\ell=-\frac{K r_s}{2} \left( u_{\ell+1}  + u_{\ell-1} \right)+c^\ast w\delta_{\ell,1}$, the equation can be written
\[
- \ell \left(\pa_\tau  u^\ast_\ell+2a u^\ast_\ell+ \frac{K r_s}{2} ( u^\ast_{\ell+1} - u^\ast_{\ell-1}) \right)=x^\ast_\ell.
\]
We have $x^\ast\in {\cal X}_{a,0}$; hence one can perform similar energy estimates to those in the proof of Lemma \ref{thm:resolvent-l1-generation} to obtain
$$ \|u^\ast \|_{a,\frac{1}{2}} \leq C \, \|x^\ast \|_{a,0}, $$
for some $C\in\R^+$. This inequality implies that $x^\ast \in {\cal X}_{a,\frac{1}{2}}$, and therefore
$$ \|u^\ast \|_{a,1} \leq C' \, \|x^\ast \|_{a,\frac{1}{2}}, $$
and \eqref{regularityu*} follows using a bootstrap argument.

Back to the main purpose of this Section, let $P_s=\text{Id}-P_0$. We claim that any $\widehat{f}$ in a sufficiently small neighborhood of the circle $\{\widehat{R}_\Theta \widehat{f_\text{s}}\}_{\Theta\in\T}$ in ${\cal X}_{a,0}$ can be written
\begin{equation} \label{ansatz}
\widehat{f}=\widehat{R}_\Theta\left(\widehat{f_\text{s}}+u\right),
\end{equation}
where $(\Theta,u)\in \T\times P_s({\cal X}_{a,0})$ is such that $\|u\|_{a,0}\to 0$ when the distance $d\left(\widehat{f},\{\widehat{R}_\Theta \widehat{f_\text{s}}\}_{\Theta\in\T}\right)\to 0$.

To see this, consider the map $F:\T\times {\cal X}_{a,0}\to\R$ defined by
$$ F(\Theta,\widehat{f}) = \text{Re}  \langle \widehat{R}_{-\Theta} \widehat{f}  - \widehat{f_\text{s}} ,  u^\ast \rangle_{a,0}, $$
which is such that $F(\Theta,\widehat{f}) =0$ iff $\widehat{f}$ satisfies \eqref{ansatz}.

We compute
$$ F(0,\widehat{f_\text{s}}) = 0\ \text{and}\ \pa_\Theta F(0, \widehat{f_\text{s}}) = - \text{Re} \left(\langle  D\widehat{R}\widehat{f_\text{s}} , u^\ast \rangle_{a,0}\right)  \neq 0, $$
hence by the Implicit Function Theorem, for $\widehat{f}$ close enough to  $\widehat{f_\text{s}}$, there is a smooth $\Theta_0 = \Theta_0(\widehat{f})$ near $0$ such that $F(\Theta_0(\widehat{f}), \widehat{f}) = 0$. This proves the claim in the neighborhood of  $\widehat{f_\text{s}}$, with $u=\widehat{R}_{-\Theta_0(\widehat{f})} \widehat{f}  - \widehat{f_\text{s}}$. To extend that property to a neighborhood of $\{\widehat{R}_\Theta \widehat{f_\text{s}}\}_{\Theta\in\T}$, notice that letting
\[
\Theta_\text{min}= \text{Arg}\min_{\Theta\in\T}\| \widehat{f} - \widehat{R}_\Theta \widehat{f_\text{s}} \|_{a,0},
\]
the element $\widehat{R}_{-\Theta_\text{min}} \widehat{f}$ is close to $\widehat{f_\text{s}}$ when $\widehat{f}$ is close to the circle. Hence, one can apply the previous argument to $\widehat{R}_{-\Theta_\text{min}} \widehat{f}$ to obtain $\widehat{f}=\widehat{R}_{\Theta_0(\widehat{R}_{-\Theta_\text{min}}\widehat{f})+\Theta_\text{min}}\left(\widehat{f_\text{s}}+u\right)$ with $u$ as desired.

\subsection{Proof of Theorem \ref{thmNL}}
By inserting the expression \eqref{ansatz} into the equation \eqref{KF} of the Kuramoto dynamics in Fourier variables, one gets after using equivariance
\[
D\widehat{R}(\widehat{f_\text{s}}+u)\frac{d\Theta}{dt}+\partial_t u=Lu+Qu.
\]
Applying $P_0$ and $P_s$ respectively, and using $P_0L=0$, $P_0L=LP_0$ and the normalization $\langle D\widehat{R}\widehat{f_\text{s}},u^\ast\rangle=\tfrac12$, two independent equations result for the variables $\Theta$ and $u$, namely
\begin{equation} \label{EQTHETA}
\frac{d\Theta}{dt}= \frac{2\text{Re}\langle Qu,u^\ast \rangle_{a,0}}{1+2\text{Re} \langle D\widehat{R} u,u^\ast \rangle_{a,0}},
\end{equation}
and
\begin{equation}
\partial_t u=Lu+ P_s Q' u\ \text{where}\ Q' u = Q u -  \frac{ 2\text{Re} \langle Qu,u^\ast \rangle_{a,0}}{1+2\text{Re} \langle D\widehat{R} u,u^\ast \rangle_{a,0}} D\widehat{R} u.
 \label{NLEQPROJECT}
\end{equation}
As intended, the right hand sides of these equations do not depend on the angular variable $\Theta$. Moreover, the Cauchy-Schwarz inequality implies
\[
\left|\langle D\widehat{R} u , u^\ast \rangle_{a,0}\right| \leq 2 \|u\|_{a,0} \| u^\ast \|_{a,1},
\]
and the property \eqref{regularityu*} implies $\| u^\ast \|_{a,1} <+\infty$. Therefore, these equations are well-defined as long as $\|u\|_{a,0}$ is small enough (so that the denominators do not vanish).

Now, if the (restriction to $\N\times\R$ of the) Fourier transform $\widehat{f_\text{in}}$ of an initial probability measure $f_\text{in}$ is sufficiently close to $\widehat{f_\text{s}}$, then not only the corresponding initial $u_\text{in}\in P_s({\cal X}_{a,0})$ is small, but the solution $\widehat{f(t)}$ must remain close to $\widehat{f_\text{s}}$ for $t\in (0,T)$, a sufficiently small time interval, by the continuous dependence in time (Proposition \ref{WP}). Hence, both ansatz \eqref{ansatz} holds and the equations above are well-defined over $(0,T)$. That these properties holds for all times (provided that $\widehat{f_\text{in}}$ is taken even closer to $\widehat{f_\text{s}}$) is a direct consequence of the following statement.

\begin{Pro}
Under condition \eqref{STABCOND}, there exist $\epsilon',b,C>0$ such that for  all $u_{\text{\rm in}} \in P_s ({\cal X}_{a,0})$ satisfying $\| u_{\text{\rm in}} \|_{a,0} < \epsilon'$,  equation \eqref{NLEQPROJECT} has a unique solution $t\mapsto u(t)$ satisfying $u(0)=u_{\text{\rm in}}$ and
\[
\|u(t)\|_{a,0}\leq C\|u_{\text{\rm in}}\|_{a,0} \, e^{-b t},\ \forall t\in\R^+.
\]
\label{ASYMPPROJECT}
\end{Pro}
This statement is not as obvious as it may look because the quadratic term $Q'$ maps ${\cal X}_{a,0}$ into ${\cal X}_{a,-1}$. The proof is given in Section \ref{S-ASYMPPROJECT} below.

In addition to ensuring that both ansatz \eqref{ansatz} and equations \eqref{EQTHETA} and \eqref{NLEQPROJECT} are globally well-defined when starting sufficiently close to $\widehat{f_\text{s}}$, Proposition \ref{ASYMPPROJECT} implies that the solution must asymptotically approach the PLS circle.
To complete the proof of Theorem \ref{thmNL}, it remains to show that the solution's angle asymptotically converges $\Theta(t)$.
We have
\[
\left|\langle Qu , u^\ast\rangle_{a,0}\right|\leq \|Qu\|_{a,-1}\| u^\ast \|_{a,1},
\]
and the definition of $Q$ and the Sobolev bound \eqref{SOBEMB}  yield
\begin{equation}
\|Qu\|_{a,-1}\leq C'  \|  u\|_{a,-\frac{1}{2}} \|u\|_{a,0}  \leq C'   \|u\|_{a,0}^2,
\label{BOUNDQ}
\end{equation}
for some $C'\in\R^+$.
Hence the driving term in equation \eqref{EQTHETA}  must also decay exponentially with rate $b$. Consequently, the following limit exists
\[
\Theta_\infty:=\lim_{t\to +\infty}\Theta(t)=\Theta(0)+\int_{\R^+}\frac{2\text{Re} \langle Qu(s),u^\ast \rangle_{a,0}}{1+2\text{Re} \langle D\widehat{R} u(s),u^\ast \rangle_{a,0}} ds,
\]
and we have
\[
\|\widehat{R}_{\Theta(t)}\left(\widehat{f_\text{s}}+u(t)\right)-\widehat{R}_{\Theta_\infty}f_\text{s}\|_{a,0}\leq |\Theta(t)-\Theta_\infty|\|\widehat{f_\text{s}}\|_{a,\frac12}+\|u(t)\|_{a,0}=O(e^{-bt}),
\]
as desired.

\subsection{Analysis of the forced linear equation}
Proposition \ref{ASYMPPROJECT} will follow from a similar result for the corresponding forced linear equation.
As suggested above, the crucial point is to show, via a suitable adaptation of the Gearhart-Pr\"uss Theorem, that the corresponding semigroup improves regularity.

To see this, given an Hilbert space $H$ with norm $\|\cdot\|_H$, a positive real number $\gamma$, and a mapping $w:\R\to H$, consider the norm $\|w\|_{H,\gamma}$ defined by
\[
\|w\|_{H,\gamma}=\left(\int_{\R^+}e^{2\gamma t}\|w(t)\|_{H}^2dt\right)^{\frac12}.
\]

We have the following statement.
\begin{Lem}
Let $X$ and $Y$ be Hilbert spaces, where $X$ is continuously embedded in $Y$. Let $A$ be a densely defined linear operator that generates a semigroup, both on $X$ and on $Y$. Assume the existence of $\gamma\in\R^+$ such that the resolvent of $A$ over both
  spaces contains the half-plane $\text{Re}(\lambda) \geq -\gamma$ and satisfies
  \begin{equation*}
    \sup_{y \in \R}  \| ((-\gamma + i y) \text{\rm Id} - A)^{-1}
    \|_{Y \to X} <C_R,
  \end{equation*}
for some $C_R\in\R^+$. Then the unique mild solution $w \in C(\R^+,Y)$ of the initial value problem
\[
\frac{dw}{dt}=Aw+G
\]
where the forcing $G : \R^+ \mapsto Y$ satisfies $\|G\|_{Y,\gamma}<+\infty$ and the initial condition $w(t)=w_\text{\rm in}$ satisfies $\|w_\text{\rm in}\|_X<+\infty$, has the following properties
\begin{itemize}
  \item[$\bullet$] $w(t) \in X$ for a.e.\ $t \in \R^+$
  \item[$\bullet$] $\|w\|_{X,\gamma}\leq C\left(\|w_\text{in}\|_X +\|G\|_{Y,\gamma}\right)$
  \end{itemize}
for some $C\in\R^+$.
\end{Lem}
{\sl Proof.} The mild solution of the initial value problem is characterized by the Duhamel's formula
\begin{equation*}
w(t) = e^{tA} w_\text{\rm in} + I(t),\ \forall t\in\R^+,\ \text{where}\ I(t)=\int_0^t e^{(t-s)A} G(s) d s.
\end{equation*}
By the Gearhart-Pr\"uss Theorem, the resolvent estimate shows that
there exists $\gamma_G > \gamma$ and $C_G\in\R^+$ such that
\begin{equation*}
  \| e^{tA} w\|_{X} \leq C_G e^{-2\gamma_G t}\|w\|_{X},\ \forall w\in X,t\in\R^+,
\end{equation*}
which yields
\begin{equation*}
  \| e^{tA} w_{\text{\rm in}} \|_{X,\gamma} \leq \frac{C_G}{2(\gamma_G-\gamma)} \| w_{\text{\rm in}}\|_X.
\end{equation*}
Moreover, for $\text{Re}(z) > - \gamma$, the Laplace transform
$\text{Lap } I$ of the integral term $I$ exists as Bochner integral over $Y$ and
satisfies
\begin{equation*}
  (\text{Lap } I)(z) = (z\text{\rm Id}-A)^{-1} (\text{Lap }G)(z),
\end{equation*}
where $\text{Lap } G$ is the Laplace transform of $G$.  On the line $\text{Re}(z) = - \gamma$, $\text{Lap } G$ exists as a $L^2$
function by the Plancherel's Theorem, and we have
\begin{equation*}
\int_\R \| (\text{Lap } G)(-\gamma+iy) \|^2_Y dy
  \leq 2\pi \| G \|^2_{Y,\gamma}.
\end{equation*}
The assumption on the resolvent estimate then implies
\begin{equation*}
\int_\R \| (\text{Lap } I)(-\gamma+iy) \|^2_X dy
  \leq 2\pi C_R\| G \|^2_{Y,\gamma}.
\end{equation*}
Using the Plancherel's Theorem again, this time to $(\text{Lap } I)(-\gamma+i\cdot)$, it follows that $I(t)
\in X$ for a.e. $t \in \R^+$ and $\| I \|^2_{X,\gamma} \leq 2\pi C_R \| G
\|^2_{Y,\gamma}$. Combined with the estimate on the initial term this shows
the claimed result. \hfill $\Box$

Now, Proposition \ref{thmL} implies that when the stability condition \eqref{STABCOND} holds, the operator ${\cal L}$ satisfies the condition of the Lemma with $X={\cal P}_s({\cal X}_{a,0}^2)$, $Y={\cal P}_s({\cal X}_{a,-1}^2)$ and $\gamma=b\in (0,a)$. The Lemma then yields the following conclusion for the initial value problem
\begin{equation}
\partial_t u=Lu+ P_sF(t),\ \text{and}\ u(0)=u_\text{\rm in}\in P_s ({\cal X}_{a,0}),
\label{FORCED}
\end{equation}
where we use the notation
\[
\|u\|_{a,k,b}=\left(\int_{\R^+}e^{2bt}\|u(t)\|_{a,k}^2dt\right)^{\frac12}.
\]
\begin{Cor}
Under the stability condition \eqref{STABCOND}, there exist $b,C>0$ such that, for every forcing signal satisfying $\|F\|_{a,-1,b}<+\infty$,
the initial value problem \eqref{FORCED} has a unique mild solution $t\mapsto u(t)\in C(\R^+ , {\cal X}_{a,-1})$ with the following properties
\begin{itemize}
\item[$\bullet$] $u(t)\in {\cal X}_{a,0}$ for a.e.\ $t\in\R^+$,
\item[$\bullet$] $\|u\|_{a,0,b}\leq C\left(\|u_\text{in}\|_{a,0}+\|F\|_{a,-1,b}\right)$.
\end{itemize}
\label{CONTROLFORCING}
\end{Cor}
%

\subsection{Proof of Proposition \ref{ASYMPPROJECT}}\label{S-ASYMPPROJECT}
The proof proceeds through a  localization of the nonlinearity.  We first show that the quantity $\|u\|_{a,0,b}$ can be made arbitrarily small for the localized system. Then we obtain a similar bound for a $L^\infty$-norm in time, which allows us to get rid of the localization and conclude the same control for the original problem.

In more details, given $\epsilon>0$, let $Q'_\epsilon:{\cal X}_{a,0}\to {\cal X}_{a,-1}$ be a smooth mapping such that
\[
Q'_\epsilon u=\left\{\begin{array}{ccl}
Q' u&\text{if}&\|u\|_{a,0}\leq \epsilon\\
0&\text{if}&\|u\|_{a,0}\geq 2\epsilon.
\end{array}\right.
\]
When $\epsilon$ is small enough, the denominator in the expression of $Q'$ remains positive over the ball $\{u\ :\ \|u\|_{a,0}< 2\epsilon\}$; hence $Q'_\eps$ is globally defined over ${\cal X}_{a,0}$ in this case.
Moreover,  using the inequality \eqref{BOUNDQ}, we infer
\[
\|Q'_\epsilon u\|_{a,-1}\leq 2\epsilon C_K  \|u\|_{a,0}.
\]
Adapting the analysis in Section \ref{S-WELLPOSED}, one can show for all $u_\text{in} \in P_s ({\cal X}_{a,0})$, there exists a unique global in time weak solution
$$ u \in C([0,T], {\cal X}_{a,0}) \cap L^2(0,T, {\cal X}_{a,\frac{1}{2}}) \quad \forall T  >0,$$
of
\begin{equation} \label{NLEQPROJECTeps}
\pa_t u =  L u + P_s Q'_\eps u
\end{equation}
with $u(0) = u_\text{in}$.  Moreover, by standard arguments, it coincides with the mild solution of equation \eqref{FORCED} with $F = Q'_\eps$.

 Applying Lemma \ref{CONTROLFORCING}, we conclude that the solution of \eqref{NLEQPROJECTeps}
  satisfies the inequality
\begin{equation}
\|u\|_{a,0,b}\leq \frac{C}{1-2\epsilon C C_K}\|u_\text{in}\|_{a,0},
\label{L2BOUND}
\end{equation}
provided that $\epsilon$ is small enough, so that the denominator here is positive.

In order to get an $L^\infty$ bound, we directly perform an estimate on equation \eqref{NLEQPROJECTeps}. We get
\[
\frac{1}{2}\frac{d}{dt}\|u \|_{a,0}^2 +
2a\|u\|_{a,\frac{1}{2}}^2\leq C_1 \|u\|_{a,0}^2 + C_1' \text{Re}
\langle Q'_\epsilon u , u \rangle_{a,0}
\]
for constants $C_1$ and $C_1'$. The second term in the right hand side
can controled as follows
\begin{equation*}
|\langle Q' u , u \rangle_{a,0}| \leq C_2 \chi\left(\frac{\|u\|_{a,0}}{\epsilon}\right) (|u(1,0)| + \| u \|_{a,0}) \| u \|_{a,\frac{1}{2}}^2
\end{equation*}
for some $C_2\in\R^+$ and where $\chi:\R^+\to [0,1]$ is a smooth function such that
\[
\chi(x)=\left\{\begin{array}{ccl}
1&\text{if}&x\leq 1\\
0&\text{if}&x\geq 2.
\end{array}\right.
\]
For $\eps$ small enough, this term can be absorbed by the left-hand side, and the following inequality results
$$ \frac{1}{2}\frac{d}{dt}\|u \|_{a,0}^2 + \frac{a}{2}\|u\|_{a,\frac{1}{2}}^2\leq C_3 \|u\|_{a,0}^2  $$
for some $C_3\in\R^+$ and then
$$  \frac{d}{dt} (e^{2bt}\|u \|_{a,0}^2)  \leq (2b + C_3) e^{2bt}  \|u\|_{a,0}^2.  $$
The Gronwall's Lemma yields in turn
$$ e^{2bt}\|u(t) \|_{a,0}^2 \leq  \|u_\text{in}\|^2_{a,0} \exp\left((2b+C_3) \int_0^t e^{2bs} \| u(s) \|^2_{a,0} ds\right). $$
Using the bound \eqref{L2BOUND}, the desired exponential decay follows
\[
\sup_{t\in\R^+}e^{2bt}\|u(t)\|_{a,0}\leq C_4 \|u_\text{in}\|^2_{a,0},
\]
for the solution of equation \eqref{NLEQPROJECTeps}. Finally, by choosing $\|u_\text{in}\|_{a,0}$ small enough, this inequality implies in particular that $\|u(t)\|_{a,0}\leq\epsilon$ for all $t\in\R^+$ and hence we have $Q_\epsilon u(t)=Q u(t)$ for all $t$, {\sl i.e.} $t\mapsto u(t)$ is actually a solution of \eqref{NLEQPROJECT}. The proof of Proposition \ref{ASYMPPROJECT} is complete.

\section{Stability condition: analysis and examples}\label{S-STAB}
As shown in Section \ref{S-LINEAR}, the stability criterion \eqref{STABCOND} in Theorem \ref{thmNL} is equivalent to the linear stability of the circle $\{{R}_\Theta {f_\text{s}}\}_{\Theta\in\T}$, more precisely that
0 is the only eigenvalue, which is simple, in the half-plane $\text{Re}(\lambda)\geq 0$, and that the rest of the spectrum lies in $\text{Re}(\lambda)\leq -\epsilon$ for some $\epsilon>0$.
Of note, that 0 must always be an eigenvalue is a consequence of the rotation symmetry $R_\Theta$ (as shown in relation \eqref{0EIGEN}). However, this property can be obtained independently, as in \cite{OW13}, by using the equations
\begin{equation}
\beta^2+2iz\beta-1=0,
\label{EQBETA}
\end{equation}
and $K\int_{\R}g(Kr_\text{s}\omega)\beta(\omega)d\omega=1$ (the latter is a rewriting using \eqref{BETA} of the self-consistency condition of the PLS $f_\text{s}$ in Section \ref{S-PLS}).
Indeed, one directly checks that
\[
\frac{K}2\left(J_0(\lambda,r_\text{s})+2\lambda J_1(\lambda,r_\text{s})+J_2(\lambda,r_\text{s})\right)=1,\ \forall \lambda\ :\ \text{Re}(\lambda)>0,
\]
from where $\det (\text{\rm Id}-\frac{K}2 M(0,r_\text{s}))=0$ immediately follows when taking the limit $\lambda\to 0$ in $\R$.

\subsection{Symmetric frequency distributions}
PLS stability depends on context and, as for existence, various situations can occur depending on the bifurcation that generates these states. For instance, when $g$ is an even function, we have $\overline{J_k(\bar\lambda,r)}=J_k(\lambda,r)$ for all $\lambda\in\C$ and then
\[
\det \left(\text{\rm Id}- \frac{K}2 M(\lambda,r_\text{s})\right)=\left(1-\frac{K}2 \left(J_0(\lambda,r_\text{s})- J_2(\lambda,r_\text{s})\right)\right)\left(1-\frac{K}2 \left(J_0(\lambda,r_\text{s})+J_2(\lambda,r_\text{s})\right)\right).
\]
Moreover, one can  show (we skip the tedious computation for brevity) that
\[
J_0(\lambda,r_\text{s})-J_2(\lambda,r_\text{s})=2h_c(\lambda) \quad\text{and}\quad
J_0(\lambda,r_\text{s})+ J_2(\lambda,r_\text{s})=2h_s(\lambda),
\]
where the functions $h_c$ and $h_s$ are defined in \cite{MS07}. In this way, we can link our stability criterion to the results of \cite{MS07}. In the case of unimodal $g$, Proposition 4 in this paper implies that for $K > K_c := \frac{2}{\pi g(0)}$, $1 - K h_c$ does not vanish over $\text{Re} (\lambda) \geq 0$, while the only zero of $1 - K h_s$ in $\text{Re}(\lambda) \geq 0$ is $\lambda = 0$. It follows that
\[
\det \left(\text{\rm Id}-\frac{K}2 M(\lambda,r_\text{s})\right)>0,\ \forall \lambda\neq 0\ \text{with}\ \text{\rm Re}(\lambda)\geq 0,
\]
for the unique PLS $f_\text{s}$ which exists for $K>K_c$. To check the second point in \eqref{STABCOND},
we use the expression of $h_s$ given in \cite{MS07}, we find
\begin{align*}
h_s'(0) & = \int_{|\omega| \geq K r_\text{s}} \frac{g(\omega)d\omega}{\sqrt{\omega^2 - (K r_\text{s})^2} (|\omega| + \sqrt{\omega^2 - (K r_\text{s})^2})} \: - \: \frac{1}{(Kr_\text{s})^2} \int_{|\omega| \leq K r_\text{s}} g(\omega) d\omega \\
& = \frac{2}{K r_\text{s}} \left(  \int_1^{+\infty}   \frac{g( K r_\text{s} \xi)d\xi}{\sqrt{\xi^2 - 1} (\xi + \sqrt{\xi^2 - 1})} - \int_0^1 g( K r_\text{s} \xi) d\xi \right) \\
& > \frac{2}{K r_\text{s}} g(K r_\text{s}) \left( \int_1^{+\infty}   \frac{d\xi}{\sqrt{\xi^2 - 1} (\xi + \sqrt{\xi^2 - 1})} - \int_0^1  d\xi \right)
\end{align*}
the last inequality coming from the fact that $g$ is unimodal. A simple computation shows that
$$ \int_1^{+\infty}   \frac{d\xi}{\sqrt{\xi^2 - 1} (\xi + \sqrt{\xi^2 - 1})}  = 1$$
so that $h'_s(0) \neq 0$.  Theorem \ref{thmNL} implies that, when it exists, this stationary solution is always asymptotically stable.

Finally, notice that uniqueness of a PLS circle does not necessarily imply its stability. Counter-examples exist.

\subsection{Stability in the Ott-Antonsen manifold}\label{S-OA}
Our next example is when $g$ is the bi-Cauchy distribution. Prior to presenting this case, we consider the dynamics in a remarkable invariant subset, the so-called Ott-Antonsen (OA) manifold \cite{OA08}, defined as the set of probability measures $f$ whose Fourier coefficients associated with the angle variable, and defined by
\begin{equation}
\widetilde{f}_\ell(\omega)=\int_\T e^{-i\ell \theta}f(d\theta,\omega),
\label{FOURIERCOEF}
\end{equation}
write
\[
\widetilde{f}_\ell(\omega)=g(\omega)h^\ell(\omega), \forall (\ell,\omega)\in \N\times\R.
\]
This set is invariant under the Kumaroto flow and the amplitude $h$ evolves according to the equation
\begin{equation}
\partial_t h(\omega)+i\omega h(\omega)+\frac{K}{2}\left(re^{-i\Omega}h^2(\omega)-re^{i\Omega}\right)=0,
\label{OADYNAM}
\end{equation}
where the order parameter $re^{i\Omega}$ can be defined in this context as
\[
re^{i\Omega}=\int_\R\overline{h(\omega)}g(\omega)d\omega.
\]
If, in addition, $g$ is meromorphic in the lower half-plane $\text{Im}(z) < 0$ with finitely many poles and sufficient decay, and if the amplitude $h$ is analytic in the same region, then the OA manifold dynamics is effectively governed by a finite dimensional system of coupled ODEs.

Now, as shown in expression \eqref{FOURCOEF} in Appendix \ref{A-STAB}, the PLS $f_\text{s}$ belongs to the OA manifold and its amplitude is given by $h_\text{s}(\omega)=\beta(\frac{\omega}{Kr_\text{s}})$ (NB: Imposing that the amplitude $h$ be analytic in the domain $\text{Im}(z) < 0$, $f_\text{s}$ turns out to be the only PLS that belongs to the OA manifold.) As before, stationary solutions of the OA dynamics \eqref{OADYNAM} come in the form of circles $\{e^{i\Theta} \beta(\frac{\omega}{Kr_\text{s}})\}_{\Theta\in\T}$.

As stability is concerned, we observe that the matrix $M(\lambda,r)$ in Theorem \ref{thmNL} satisfies the conjugacy equation
\[
\frac{K}2 M(\lambda,r)=PB(\lambda)P^{-1}\quad \text{where}\quad P=\left(\begin{array}{cc}
1&i\\
1&-i\end{array}\right)
\]
and $B(\lambda)$ is the matrix in \cite{OW13}, involved in the PLS stability condition in the OA manifold. (Of note, the analysis in \cite{OW13} deals with strong topology on amplitude functions $h$ and faces the issue of a continuous spectrum on the imaginary axis. Yet, the strategy for spectral stability is similar to the one developed here.) It follows that the stability criterion in this invariant subset, namely
\[
\left\{\begin{aligned}
&\det \left(\text{\rm Id}-B(\lambda)\right)\neq 0,\ \forall \lambda \neq 0\ \text{with}\ \text{\rm Re}(\lambda)\geq 0,\\
&0\ \text{is a simple zero of}\ \lambda\mapsto \det \left(\text{\rm Id}-B(\lambda)\right),
\end{aligned}
\right.
\]
coincides with the condition \eqref{STABCOND} in the full space. In other words, no loss of generality results in investigating the existence and stability of $f_\text{s}$ in the OA manifold. The Ott-Antonsen ansatz is perfectly legitimate.

\subsection{Existence and stability of PLS for bi-Cauchy frequency distributions}
Here, we consider the existence and stability of the PLS for the bi-Cauchy frequency distribution, which has been frequently employed in illustrations of the Kuramoto dynamics. This distribution is defined by
\[
g_{\Delta,\omega_0}(\omega)=\frac{\Delta}{2\pi}\left(\frac1{(\omega-\omega_0)^2+\Delta^2}+\frac1{(\omega+\omega_0)^2+\Delta^2}\right),
\]
where $\omega_0,\Delta\in\R^+$ and we take $\omega_0>\frac{\Delta}{\sqrt{3}}$ (so that the distribution is bimodal, not unimodal). The existence statement relies on the convex map $\Psi$ defined by
\[
\Psi_{\Delta,\omega_0}(x)=2\Delta\left(\frac1{1-x}+\frac{\omega_0^2}{\Delta^2}\frac{1-x}{1+x}^2\right),\ \text{for}\ x\in [0,1).
\]
Notice that $\Psi_{\Delta,\omega_0}(0)=\frac{\pi}{2g_{\Delta,\omega_0}(0)}$ and $\Psi_{\Delta,\omega_0}'(0)<0$ for the derivative. The analysis of the existence and stability condition in the OA manifold yields the following conclusion.
\begin{Exam}
Let $\Delta>0$ and $\omega_0>\frac{\Delta}{\sqrt{3}}$ be fixed. Then, the existence and stability of the circle $\{R_\Theta f_\text{s}\}_{\Theta\in\T}$ can be enumerated as follows, depending on $K$:
\begin{itemize}
\item[$\bullet$] If $K<{\displaystyle\min\limits_{x\in [0,1)}}\Psi_{\Delta,\omega_0}(x)$, then no PLS exist.
\item[$\bullet$] If $K\in \left({\displaystyle\min\limits_{x\in [0,1)}}\Psi_{\Delta,\omega_0}(x),\frac{\pi}{2g_{\Delta,\omega_0}(0)}\right)$, then two PLS circles exist with respective order parameter $r_-$ and $r_+$, where $r_\pm=\rho_\pm\sqrt{1-\left(\frac{4\omega_0}{K}\frac{1}{1+\rho_\pm^2}\right)^2}$ and $\rho_-<\rho_+$ are defined by $K=\Psi_{\Delta,\omega_0}(\rho_-^2)=\Psi_{\Delta,\omega_0}(\rho_+^2)$. The PLS circle associated with $r_+$ is asymptotically stable (in the sense of Theorem \ref{thmNL}), the one associated with $r_-$ is unstable.
\item[$\bullet$] For $K>\frac{\pi}{2g_{\Delta,\omega_0}(0)}$, a unique asymptotically stable PLS circle exists (whose order parameter is the continuation of the solution branch $r_+$).
\end{itemize}
 \end{Exam}
\begin{figure}[h]
\centerline{\includegraphics[scale=0.35]{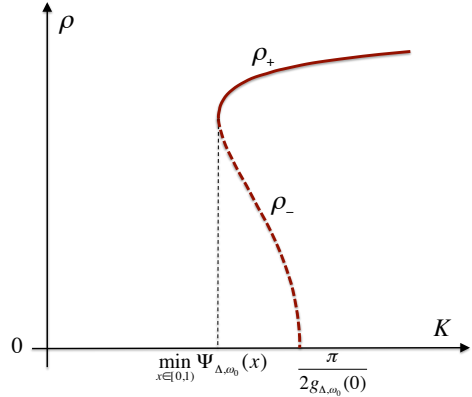}}
\caption{Schematic PLS bifurcation diagram for the bi-Cauchy distribution $g_{\Delta,\omega_0}$ with $\omega_0=2\Delta$ (obtained using the graph of the function $\Psi_{\Delta,\omega_0}$, see text).}
\label{BIFURCATIONBICAUCHY}
\end{figure}
In addition, the solution branch $K\mapsto r_+(K)$ is increasing and the other branch $K\mapsto r_-(K)$ reaches 0 for $K=\frac{\pi}{2g_{\Delta,\omega_0}(0)}$. An illustration of these branches, expressed in terms of $\rho_\pm$, is given in Figure \ref{BIFURCATIONBICAUCHY}.

{\sl Proof of existence.}  We borrow notation from \cite{MBSOSA09}.  Thanks to the OA ansatz,  it is shown there that   the existence of a PLS with order parameter $r_\text{s}\in\R^+$ is equivalent to the existence of a solution $(\rho_\text{s}, \varphi_\text{s})$ of   the equations
\begin{equation}
F_\rho(\rho_\text{s},\rho_\text{s},\psi_\text{s})=0\quad\text{and}\quad
F_\psi(\rho_\text{s},\rho_\text{s},\psi_\text{s})=0,
\label{EQUOA}
\end{equation}
where
\[
F_\rho(\rho_1,\rho_2,\psi)=-\Delta\rho_1+\frac{K}4(1-\rho_1^2)(\rho_1+\rho_2\cos\psi),
\]
and
\[
F_\psi(\rho_1,\rho_2,\psi)=2\omega_0-\frac{K}4\frac{\rho_1^2+\rho_2^2+2\rho_1^2\rho_2^2}{\rho_1\rho_2}\sin\psi.
\]
The order parameter $r_\text{s}$ is then given by $r_\text{s} = \rho_s \, |\cos(\varphi_\text{s})|$.

The equations \eqref{EQUOA} are then equivalent to $K=\Psi_{\Delta,\omega_0}(\rho_\text{s}^2)$, and  $\sin\psi_\text{s}=\frac{4\omega_0}{K}\frac{1}{1+\rho_\text{s}^2}$. The conclusion then results from the analysis of the map $\Psi_{\Delta,\omega_0}$. \hfill $\Box$

{\sl Proof of stability.} The distribution $g_{\Delta,\omega_0}$ is a rational function with two poles $\pm w_0-i\Delta$ in the lower half-plane. Hence, given any point in the OA manifold, the only contribution of corresponding amplitude $h$ to the order parameter $re^{i\Omega}$ is through the variables
\[
z_1=\overline{h(-\omega_0-i\Delta)}\quad \text{and}\quad z_2=\overline{h(\omega_0-i\Delta)}.
\]
Therefore, the OA dynamics \eqref{OADYNAM} in this case is entirely controlled by the dynamical system in $\R^4$ defined by (equations (16) and (17) in \cite{MBSOSA09})
\[
\left\{\begin{array}{l}
\dot{z}_1=-(\Delta +i\omega_0)z_1+\frac{K}4\left(z_1+z_2-(\overline{z_1}+\overline{z_2})z_1^2\right)\\
\dot{z}_2=-(\Delta -i\omega_0)z_2+\frac{K}4\left(z_1+z_2-(\overline{z_1}+\overline{z_2})z_2^2\right)
\end{array}\right.
\]
that has a circle of stationary solutions $\{e^{i\Theta} \beta(\frac{-\omega_0-i\Delta}{Kr_\text{s}}),e^{i\Theta} \beta(\frac{\omega_0-i\Delta}{Kr_\text{s}})\}_{\Theta\in\T}$, corresponding to PLS. Its stability can be investigated along the same lines as in the general case, namely by considering the linearized perturbation dynamics and the complexification of the corresponding linear operator. The analysis shows that spectral stability of this circle is given by the criterion \eqref{STABCOND} with
\[
J_k(\lambda,r)=\frac{K}4\left(\frac{\beta^k\left(\frac{\omega_0-i\Delta}{Kr}\right)}{\lambda+i\omega_0+\Delta+Kr\beta\left(\frac{\omega_0-i\Delta}{Kr}\right)}+\frac{\beta^k\left(\frac{-\omega_0-i\Delta}{Kr}\right)}{\lambda-i\omega_0+\Delta+Kr\beta\left(\frac{-\omega_0-i\Delta}{Kr}\right)}\right).
\]
Yet, spectral stability can be also checked by passing to polar coordinates in the 4-dimensional system, as in \cite{MBSOSA09}. This amounts to investigate the stability of the stationary point $(\rho_\text{s},\rho_\text{s},\psi_\text{s})$ of the flow governed by
\[
\left\{\begin{array}{l}
\dot\rho_1=F_\rho(\rho_1,\rho_2,\psi)\\
\dot\rho_2=F_\rho(\rho_2,\rho_1,\psi)\\
\dot\psi=F_\psi(\rho_1,\rho_2,\psi)
\end{array}\right.
\]
Thanks to the symmetry $(\rho_1,\rho_2,\psi)\mapsto (\rho_2,\rho_1,\psi)$, the derivative at $(\rho_\text{s},\rho_\text{s},\psi_\text{s})$ has a transverse eigenvector $(1,-1,0)$ with eigenvalue $-\Delta+\frac{K}4(1-3\rho_\text{s}^2-(1+\rho_\text{s}^2)\cos\psi_\text{s})$, which has been shown in \cite{MBSOSA09} to be negative when computed for $\rho_\text{s}=\rho_+$.

In the complement space $\text{Span}\left\{(1,1,0),(0,0,1)\right\}$, the eigenvalue equation writes $\lambda^2+2b\lambda+c=0$ where
\[
b=\Delta\frac{1+2\rho_\text{s}^2}{1-\rho_\text{s}^2}-\frac{K}{2}(1+\rho_\text{s}^2),
\]
and
\[
c=\rho_\text{s}^2\left(4\Delta^2\frac{1+\rho_\text{s}^2}{(1-\rho_\text{s}^2)^2}-4\omega_0^2\frac{1-\rho_\text{s}^2}{(1+\rho_\text{s}^2)^2}-K\Delta\frac{1+\rho_\text{s}^2}{1-\rho_\text{s}^2}\right).
\]
A sufficient condition for the roots to have negative real part is $b>0$ and $c>0$. Direct calculations show that the latter 1/ is necessary, 2/ implies the former and 3/ is equivalent to
\[
K<\Phi_{\Delta,\omega_0}(\rho_\text{s}^2),\ \text{where}\ \Phi_{\Delta,\omega_0}(x)=4\Delta\left(\frac1{1-x}-\frac{\omega_0^2}{\Delta^2}\frac{(1-x)^2}{(1+x)^3}\right).
\]
Therefore, all we need to check is $K < \Phi_{\Delta,\omega_0}(\rho_\text{s}^2)$. This inequality turns out to be equivalent to $\Psi_{\Delta,\omega_0}'(\rho_\text{s}^2)>0$ which hold for $\rho_\text{s}=\rho_+$ (and not for $\rho_\text{s}=\rho_-$). \hfill $\Box$

\section*{Acknowledgements}

We thank T.\ Gallay, G.\ Giacomin, N.\ Masmoudi and F.\ Rousset for
fruitful discussions and relevant comments. H.D.\ is supported by the
UK Engineering and Physical Sciences Research Council (EPSRC) grant
EP/H023348/1 for the University of Cambridge Centre for Doctoral
Training, the Cambridge Centre for Analysis. B.F.\ and D.G-V.\ have
received support from the CNRS PEPS program ``Physique Th\'eorique et
ses Interfaces''. D.G-V.\ also acknowledges the support of ANR project
Dyficolty ANR-13-BS01-0003-01 and the support of program
ANR-11-IDEX-005.

\appendix

\section{Partially locked states: stability requirements}\label{A-STAB}
This appendix collects preliminary properties to the stability analysis of PLS, in particular the fact that only those states with $\alpha=1\ \text{a.e.}$ can have finite $\|\cdot\|_{a,0}$-norm. Notice first that solving the stationary state equation for the Fourier coefficients defined by equation \eqref{FOURIERCOEF}, using the assumption $r_\text{pls}\in\R^+$, easily yields the following expression (which, evidently, consists of the Fourier coefficients of expression \eqref{EXPREPLS})
\begin{equation}
(\widetilde{f_\text{pls}})_\ell(\omega)=\left\{\begin{array}{ccl}
\left(\alpha(\omega)\beta^\ell\left(\frac{\omega}{Kr_\text{pls}}\right)+(1-\alpha(\omega))\beta_-^\ell\left(\frac{\omega}{Kr_\text{pls}}\right)\right)g(\omega)&\text{if}&|\omega|\leq Kr_\text{pls}\\
\beta^\ell\left(\frac{\omega}{Kr_\text{pls}}\right) g(\omega) &\text{if}&|\omega|> Kr_\text{pls}
\end{array}\right.
\label{FOURCOEF}
\end{equation}
for all $\ell\in\Z$, where $\beta$ is given in equation \eqref{BETA} and
\[
\beta_-(\omega)=-i\omega-\sqrt{1-\omega^2},\ \text{for}\ \omega\in [-1,1].
\]
Notice that $\beta_-$ is the other root of equation \eqref{EQBETA} and, thanks to the branch cut choice of $\sqrt{1-z^{-2}}$ in the unit disk, the function $\beta$ turns out to be analytic in the lower half plane $\text{Im}(z) < 0$.

Our first property is that the exponential weight on the Fourier transform  is equivalent to impose holomorphic continuation in an horizontal strip of the lower half-plane.
\begin{Lem}
Let $f$ be a complex valued Radon measure. The map $\tau\mapsto\widehat{f}(\tau)e^{a\tau}\in L^2(\R)$ for some $a>0$ iff there exists $F:\C\to\C$ such that
\begin{itemize}
\item[(i)] F  is holomorphic in the strip $\{x+iy\ :\ y\in (-a,0)\}$.
\item[(ii)] ${\displaystyle\sup\limits_{y\in [-a,-\epsilon]}} \|F(\cdot+i y)\|_{L^2(\R)} <+\infty$ and $\|F(\cdot+i y)\|_{L^2(\R)}=\|e^{-y\cdot}\widehat{f}(\cdot)\|_{L^2(\R)}$ for $y\in [-a,-\epsilon]$, for all $\epsilon>0$.
\item[(iii)] ${\displaystyle\lim\limits_{y \rightarrow 0^-}} F(\cdot + i y) = f(\cdot)$ in $\: \mathcal{S}'(\R)$.
\end{itemize}
\label{HOLOMCONT}
\end{Lem}
{\sl Proof.} Assume that $f \in \mathcal{S}'(\R, \C)$ is such that $\tau\mapsto\widehat{f}(\tau)e^{a\tau}\in L^2(\R)$ for some $a>0$ and let $F$ be defined by
\[
F(z)=\frac{1}{2\pi}\int_\R e^{iz\tau} \widehat{f}(\tau) d\tau.
\]
Item (i) is a simple consequence of holomorphy under integral sign. Item (ii) is a basic application of Plancherel isometry.

To prove (iii), we first observe that Plancherel Theorem implies the following relation
\[
\langle F(\cdot + i y), \varphi \rangle = \frac{1}{2\pi} \int_\R e^{-y \tau} \widehat{f}(\tau) \overline{\widehat{\varphi}(\tau)} d \tau,
\]
for all $\varphi \in \mathcal{S}(\R)$. A simple application of the Dominated Convergence Theorem then yields
$$
\lim_{y\to 0}\frac{1}{2\pi} \int_\R e^{-y \tau} \widehat{f}(\tau) \overline{\widehat{\varphi}(\tau)} d \tau = \frac{1}{2\pi} \int_\R  \widehat{f}(\tau) \overline{\widehat{\varphi}(\tau)} = \langle f , \varphi \rangle, $$
as desired.

Conversely, assume that (i)-(iii) are fulfilled. By (iii) and the continuity of the Fourier transform in $\mathcal{S}'$, we have
$$
\lim_{y \rightarrow 0^-} \widehat{F(\cdot + i y)} =\widehat{f},
$$
in $\mathcal{S}'$.
Now (ii) implies that we can write for almost every $\tau$
$$
\widehat{F(\cdot + i y)}(\tau) = e^{-y\tau} \int_{\R + i y} e^{-i x\tau } F(x) dx,
$$
where the integral is to be understood in semi-convergence sense. Proceeding as in \cite{R87}, holomorphy of $F$ yields
\begin{align*}
 \int_{\R + i y} e^{-i x\tau } F(x) dx  = \int_{\R - i a} e^{-i x\tau } F(x) dx = e^{-a\tau} \int_\R e^{-i x \tau} F(x- ia) dx = e^{-a\tau} u(\tau)
\end{align*}
where $u\in L^2(\R)$ as the Fourier transform of $F(\cdot- ia)\in L^2(\R)$. Therefore, we get
\begin{align*}
 \widehat{f} & = \lim_{y \rightarrow 0^-}  e^{-y \tau} e^{-a \tau} u \quad \text{in} \: \mathcal{S}'(\R) \\
& = \lim_{y \rightarrow 0^-}  e^{-y \tau} e^{-a \tau}u \quad \text{in} \: \mathcal{D}'(\R) \\
& =  e^{-a \tau} u
\end{align*}
from where the conclusion follows using $u\in L^2(\R)$.\hfill $\Box$

\begin{Pro}
Assume that $g\in L^1(\R)$ and that $\|\widehat{g}\|_{a}<+\infty$ for some $a>0$. Then, for every $r_\text{pls}\in\R^+$, the probability measure $f_\text{pls}$ defined by equation \eqref{EXPREPLS} satisfies $\|\widehat{f_\text{pls}}\|_{a,k}<+\infty$ for some $k\in\Z$, iff $\alpha=1\ \text{a.e.}$ Moreover, $\|\widehat{f_\text{pls}}\|_{a,k}<+\infty$ iff $\|\widehat{f_\text{pls}}\|_{a,k'}<+\infty$ for all $k'\neq k$.
\label{UNIQPLS}
\end{Pro}
{\sl Proof.} The proof decomposes into two steps. First, we show that $\|\widehat{f_\text{s}}\|_{a,k}<+\infty$. Then we prove that no other PLS with order parameter $r_\text{s}\in\R^+$ can have finite $\|\cdot\|_{a,k}$-norm.

{\sl Proof that $\|\widehat{f_\text{s}}\|_{a,k}<+\infty$.} From the expression of $\widehat{f_\text{s}}$ (obtained from expression \eqref{FOURCOEF} with $\alpha=1\ \text{a.e.}$), consider the quantity
\begin{equation}
\sum_{\ell\in\N}\int_\R(1+|\omega-ia|^2)|g(\omega-ia)|^2\ell^{2k}\left|\beta\left(\frac{\omega-ia}{Kr_\text{pls}}\right)\right|^{2\ell}d\omega.
\label{PLANCHEST}
\end{equation}
Equation \eqref{EQBETA} implies that $|\beta(x)|\leq 1$ and $\left|\frac1{\beta(x)}-\beta(x)\right|=2|x|$; hence given $\rho\in (0,1)$, we have $|\beta(x)|<\rho$ when $|x|>\frac{1+\rho}{2\rho}$. Using that $\beta$ is analytic in the lower half-plane, the maximum modulus principle then yields
\[
\sup_{\omega\in\R}\left|\beta\left(\frac{\omega-ia}{Kr_\text{pls}}\right)\right|<1,
\]
and the relation
\[
\int_\R(1+|\omega-ia|^2)|g(\omega-ia)|^2d\omega=\|\widehat{g}\|_{a},
\]
implies that the quantity \eqref{PLANCHEST} must be finite. However, when regarding $(\widehat{f_\text{s}})_\ell$ as a function of a real variable, Lemma \ref{HOLOMCONT} implies that each term in this sum is equal to $\ell^{2k}\|(\widehat{f_\text{s}})_\ell\|_a$. Therefore, \eqref{PLANCHEST} is nothing but $\|\widehat{f_\text{s}}\|_{a,k}$ and the conclusion immediately follows.

{\sl Proof of uniqueness.} The proof proceeds by contradiction and relies on the following statement, whose proof is given below.
\begin{Lem}
Let $f\in L^1(\R)$ be such that $f(x)=0$ for $|x|>\delta$ for some $\delta>0$, and $\tau\mapsto \widehat{f}(\tau)e^{a\tau}\in L^2(\R)$ for some $a>0$. Then $f=0$ a.e.
\label{ZERO}
\end{Lem}
Now, in addition to $f_\text{s}$, assume the existence of a PLS $f_\text{pls}$ with order parameter $r_\text{s}$, and $\alpha\neq 1$ and $\|\widehat{f_\text{pls}}\|_{a,k}<+\infty$. Then expression \eqref{FOURCOEF} implies that the first Fourier coefficient of the difference $h=f_\text{pls}-f_\text{s}$ satisfies
\begin{itemize}
\item[$\bullet$] $\widetilde{h}_1\in L^1(\R)$,
\item[$\bullet$] $\widetilde{h}_1(\omega)=0$ for $|\omega|>Kr_\text{s}$,
\item[$\bullet$] the Fourier transform satisfies $\tau\mapsto \widehat{h}_1(\tau)e^{a\tau}\in L^2(\R)$.
\end{itemize}
However, Lemma \ref{ZERO} asserts that $h=0$; hence the contradiction. \hfill $\Box$

{\sl Proof of Lemma \ref{ZERO}.}  Consider the shifted $f$ such that $\text{supp}(f) \subset [0,2\delta]$. We have $|\widehat{f}(z)| \leq \| f \|_1$ for $\text{Im} (z) \geq 0$. Now, the map $h$ be defined by
\begin{equation*}
h(z) = \widehat{f}\left(\frac{i-iz}{1+z}\right).
\end{equation*}
is holomorphic in the unit ball $D= \{z : |z| < 1\}$ and continuous in $\overline{D} \setminus \{-1\}$. Moreover, $h$ is bounded within $\overline{D} \setminus \{-1\}$ by $\|f\|_1$.

Up to dividing by $z^n$, we can assume w.l.o.g., that $h(0) \not = 0$. By contradiction, assume that $f$ is not identically 0. The mapping $z\mapsto \log |h(z)|$ is subharmonic, {\sl viz.}\ (Theorem 15.19 in \cite{R87})
\begin{equation*}
\log |h(0)| \leq \frac{1}{2\pi} \lim_{r\to 1^-}\int_\T \log|h(re^{i\theta})|d \theta=\frac{1}{2\pi} \int_\T \log|h(e^{i\theta})|d \theta,
\end{equation*}
where the equality follows from Lebesgue's Dominated Convergence Theorem based on that the quantity $|h(e^{ri\theta})|$ is bounded above. Hence the negative part $\mapsto \log_{-}|h(e^{i\cdot})|$ must be integrable over $\T$, {\sl i.e.}
\begin{equation*}
\int_\T\log_{-} |h(e^{i\theta})| d \theta = 2 \int_\R  \frac{\log_{-} |\widehat{f}(x)|}{1+x^2}d x< +\infty,
\end{equation*}
where we have used a change of variable.

Now, let $A = \{x \in \R^+: e^{ax} |f(x)| > 1\}$. We must have $\text{Leb}(A) < +\infty$, otherwise we would have $\|f\|_a=+\infty$. Moreover,
\[
\int_\R \frac{\log_{-} |f(x)|}{1+x^2}d x\geq \int_{\R^+} 1_{x \not \in A} \frac{\log_{-} |f(x)|}{1+x^2}d x
\geq \int_{\R^+} 1_{x \not \in A} \frac{ax}{2(1+x^2)} d x.
\]
Now, using that $\text{Leb}(A) < +\infty$, we get
\begin{equation*}
\int_{\R^+} 1_{x \in A} \frac{ax}{1+x^2} d x \leq C\text{Leb}(A)< +\infty,
\end{equation*}
for some $C\in\R^+$.
However, the integral $\int_{\R^+} \frac{ax}{1+x^2} d x$ diverges. Therefore, the integral $\int_{\R^+} 1_{x \not \in A} \frac{ax}{2(1+x^2)} d x$ also diverges, and this contradicts the fact that $\theta\mapsto \log_{-}|h(e^{i\theta})|$ is integrable. \hfill $\Box$

\section{Weak convergence of measures induced by $\|\cdot\|_{a,0}$}\label{A-WEAKCONV}
\begin{Lem}
Let $\{f_n\}_{n\in\N}$ be a sequence of probability measures on the cylinder with frequency marginal $g$ and let $f$ be with the same property. We have
\[
\lim_{n\to\infty}\|\widehat{f_n}-\widehat{f}\|_{a,0}=0\ \Longrightarrow\ \lim_{n\to\infty}f_n=f,
\]
where convergence here is understood in the weak sense.
\label{WEAKCONV}
\end{Lem}
{\sl Proof.} The sequence of frequency marginals $\int_\T f_n(d\theta,d\omega)$ is  tight, because it is constant. Hence, the sequence $\{f_n\}_{n\in\N}$ itself is  tight.

Let $f'$ be any accumulation point of $\{f_n\}_{n\in\N}$ and let $\{n_i\}$ be the corresponding subsequence. Convergence in weak topology implies
\[
\lim_{i\to\infty}(\widehat{f_{n_i}})_\ell(\tau)=\widehat{f'}_\ell(\tau),\ \forall (\ell,\tau).
\]
However, the convergence $\|\widehat{f_n}-\widehat{f}\|_{a,0}\to 0$ implies that every $(\widehat{f_n})_\ell$ converges in $H^1([-m,m])$, for every $m\in\R^+$. By the Sobolev embedding $H^1([-m,m])\hookrightarrow C_0([-m,m])$, this implies
\[
\lim_{n\to\infty}(\widehat{f_{n}})_\ell(\tau)=\widehat{f}_\ell(\tau),\ \forall (\ell,\tau).
\]
Since the Fourier transform is one-to-one, we must have $f'=f$ for every accumulation point $f'$. Hence $\lim_{n\to\infty}f_n=f$. \hfill $\Box$

\end{document}